\documentclass[a4paper,
fontsize=11pt,%
oneside,%
numbers=enddot]{scrartcl}
\KOMAoptions{DIV=12} 
\usepackage[usenames,dvipsnames]{pstricks}
\usepackage{pstricks-add}
\usepackage{pst-plot}
\usepackage{pst-math}
\usepackage[utf8]{inputenc}
\usepackage[T1]{fontenc}
\usepackage{textcomp}
\usepackage[fleqn]{amsmath}
\usepackage{amsthm}
\usepackage{amssymb}
\usepackage{amsfonts}
\usepackage{graphicx}
\usepackage{libertine}
\usepackage[libertine,
slantedGreek,
nosymbolsc,
nonewtxmathopt,
subscriptcorrection]{newtxmath}
\usepackage[scaled=0.95,varqu,varl]{inconsolata}
\usepackage[scr=boondox]{mathalpha} 
\usepackage{euscript} 
\usepackage{leftindex}
\allowdisplaybreaks[1] 
\numberwithin{equation}{section} 
\usepackage{xcolor}
\definecolor{dblue}{HTML}{0171BD}
\definecolor{dgreen}{HTML}{02724A}
\definecolor{dgreen2}{HTML}{025951}
\definecolor{dred}{HTML}{D90404}
\definecolor{dbrown}{HTML}{780080}
\definecolor{dviolet}{HTML}{42208C}
\definecolor{labelkey}{HTML}{025951}
\definecolor{refkey}{HTML}{025951}
\definecolor{orng}{HTML}{D35400}
\definecolor{pblue}{rgb}{0.1176,0.5647,1}
\definecolor{pgreen}{rgb}{0.1961,0.8039,0.1961}
\definecolor{pred}{rgb}{1.0,0.2706,0.0}
\definecolor{fred}{rgb}{0.98,0.40,0.93}
\definecolor{pyellow}{rgb}{1.0,0.6471,0.0}
\definecolor{color96b}{rgb}{0.9,0.90,1.0}
\definecolor{red1}{rgb}{0.0,0.60,0.6}
\definecolor{red2}{rgb}{0.90,0.0,0.0}
\usepackage{tikz}
\usepackage{pgfplots}
\usepgfplotslibrary{colormaps}

\pgfplotsset{colormap={setting1}{color=(pblue) color=(pyellow) 
color=(pred)},
colormap={setting2}{rgb255=(255,0,0) rgb255=(255,255,0)}}

\addtokomafont{section}{\centering}

\usepackage{enumitem}
\setlist{itemsep=-2.0pt}

\usepackage{upref} 
\usepackage{hyperref}
\hypersetup{colorlinks=true,
linktocpage=true,
linkcolor=dblue,
citecolor=dgreen,
urlcolor=dred,
pdfencoding=auto,
hypertexnames=false}
\makeatletter
\g@addto@macro\th@plain{
\thm@headfont{\bfseries\sffamily}
\thm@notefont{}}
\g@addto@macro\th@definition{
\thm@headfont{\bfseries\sffamily}
\thm@notefont{}}
\g@addto@macro\th@remark{
\thm@headfont{\bfseries\sffamily}
\thm@notefont{}}
\makeatother
\theoremstyle{plain}
\newtheorem{theorem}{Theorem}[section]
\newtheorem{proposition}[theorem]{Proposition}

\newtheorem{lemma}[theorem]{Lemma}

\theoremstyle{definition}
\newtheorem{definition}[theorem]{Definition}

\newtheorem{problem}[theorem]{Problem}
\newtheorem{assumption}[theorem]{Assumption}
\theoremstyle{remark}
\newtheorem{remark}[theorem]{Remark}

\usepackage[extdef=true]{delimset}
\DeclareMathDelimiterSet{\scal}[2]{
\selectdelim[l]<{#1}
\mathpunct{}\selectdelim[p]|
{#2}\selectdelim[r]>}

\newcommand{\menge}[2]{\left\{{#1}\mid{#2}\right\}} 
\DeclareMathDelimiterSet{\Menge}[2]{\selectdelim[l]\{
{#1}\selectdelim[m]|{#2}\selectdelim[r]\}}
\newcommand*\Cdot{{\mkern 1.6mu\cdot\mkern 1.6mu}}
\makeatletter
\def\upintkern@{\mkern-7mu\mathchoice{\mkern-3.5mu}{}{}{}}
\def\upintdots@{\mathchoice{\mkern-4mu\@cdots\mkern-4mu}%
{{\cdotp}\mkern1.5mu{\cdotp}\mkern1.5mu{\cdotp}}%
{{\cdotp}\mkern1mu{\cdotp}\mkern1mu{\cdotp}}%
{{\cdotp}\mkern1mu{\cdotp}\mkern1mu{\cdotp}}}
\makeatother
\DeclareFontFamily{OMX}{mdbch}{}
\DeclareFontShape{OMX}{mdbch}{m}{n}{ <->s * [0.8] mdbchr7v }{}
\DeclareFontShape{OMX}{mdbch}{b}{n}{ <->s * [0.8] mdbchb7v }{}
\DeclareFontShape{OMX}{mdbch}{bx}{n}{<->ssub * mdbch/b/n}{}
\DeclareSymbolFont{uplargesymbols}{OMX}{mdbch}{m}{n}
\SetSymbolFont{uplargesymbols}{bold}{OMX}{mdbch}{b}{n}
\DeclareMathSymbol{\upintop}{\mathop}{uplargesymbols}{82}
\DeclareMathSymbol{\upointop}{\mathop}{uplargesymbols}{"48}
\makeatletter
\renewcommand{\int}{\DOTSI\upintop\ilimits@}
\renewcommand{\oint}{\DOTSI\upointop\ilimits@}
\makeatother

\newcommand{\qq}{\mathscr{Q}}
\newcommand{\RR}{\mathbb{R}}
\newcommand{\NN}{\mathbb{N}}
\newcommand{\HH}{\mathcal{H}}
\newcommand{\GG}{\mathcal{G}}

\newcommand{\FF}{\EuScript{F}}

\newcommand{\GW}{\mathsf{G}_{\omega}}

\newcommand{\gw}{\mathsf{g}_{\omega}}

\newcommand{\LW}{\mathsf{L}_{\omega}}

\newcommand{\HS}{\mathsf{H}}
\newcommand{\dft}{\mathsf{DFT}}
\newcommand{\idft}{\mathsf{IDFT}}
\newcommand{\pinf}{{+}\infty}
\newcommand{\minf}{{-}\infty}

\newcommand{\RXX}{\intv{\minf}{\pinf}}
\newcommand{\RX}{\intv[l]0{\minf}{\pinf}}
\newcommand{\RP}{\intv[r]0{0}{\pinf}}
\newcommand{\RPP}{\intv[o]0{0}{\pinf}}

\newcommand{\emp}{\varnothing}

\newcommand{\minimize}[2]{\underset{\substack{{#1}}}
{\operatorname{minimize}}\;\;#2}

\newcommand{\pushfwd}%
{\ensuremath{\mbox{\Large$\,\triangleright\,$}}}
\newcommand{\epushfwd}%
{\ensuremath{\mbox{\Large$\,\trianglerightdot\,$}}}

\DeclareMathOperator{\Argmin}{Argmin}
\DeclareMathOperator*{\argmin}{argmin}
\newcommand{\Id}{\mathrm{Id}}

\DeclareMathOperator{\ran}{range}

\DeclareMathOperator{\dom}{dom}

\DeclareMathOperator{\sign}{sign}

\DeclareMathOperator{\prox}{prox}
\DeclareMathOperator{\proj}{proj}

\newcommand{\lenv}[2]{\mathrm{lenv}_{\!#2}{\,#1}}
\newcommand{\uenv}[2]{\mathrm{uenv}_{\!#2}{\,#1}}
\DeclareMathOperator{\ave}{\mathrm{ave}}
\DeclareMathOperator{\cav}{\mathrm{cav}}
\newcommand{\pcmx}[1]{\mathrm{pcm}_{#1}}
\newcommand{\spcm}[1]{\mathfrak{m}_{#1}}
\newcommand{\pave}[1]{\mathrm{pav}_{\!#1}}
\newcommand{\scav}{\mathfrak{c}}

\newcommand{\lev}[1]%
{{\ensuremath{{{{\operatorname{lev}}}_{\leq #1}}\,}}}

\newcommand{\smalldot}{\ensuremath{\raisebox{0.12ex}{\scalebox{0.8}
{$\cdot$}}}}\def\trianglerightdot{{\mkern+1mu\smalldot\mkern-4.5mu 
\triangleright}}

\DeclareFontFamily{U}{mathb}{}
\DeclareFontShape{U}{mathb}{m}{n}{<-5.5> mathb5 <5.5-6.5> mathb6 
<6.5-7.5> mathb7 <7.5-8.5> mathb8 <8.5-9.5> mathb9 <9.5-11> mathb10
<11-> mathb12}{}
\DeclareSymbolFont{mathb}{U}{mathb}{m}{n}
\DeclareFontSubstitution{U}{mathb}{m}{n}
\DeclareMathSymbol{\blackdiamond}{\mathbin}{mathb}{"0C}

\renewcommand{\leq}{\leqslant}
\renewcommand{\geq}{\geqslant}

\newcommand{\weakly}{\rightharpoonup}

\newcommand{\proxcc}[1]{{\ensuremath{\overset{#1}{\blackdiamond}}}}

\newcommand{\Rcm}[1]{\overset{\mathord{\blackdiamond}}{\mathsf{M}}_{#1}}

\renewenvironment{abstract}{%
\vspace*{-0.50cm}
\small
\quotation%
\noindent%
{\normalfont\bfseries\sffamily
\nobreak\abstractname\ }%
}{%
\endquotation%
\medskip
}
\renewcommand{\abstractname}{Abstract.}
\newcommand\keywordsname{Keywords.}
\newenvironment{keywords}
{\renewcommand\abstractname{\keywordsname}\begin{abstract}}
{\end{abstract}}
\usepackage[auth-sc]{authblk}
\newcommand{\email}[1]{\href{mailto:#1}{\nolinkurl{#1}}}
\renewcommand*\Affilfont{\normalfont\normalsize}
\newcommand\affilcr{\protect\\ \protect\Affilfont}
\makeatletter
\renewcommand\AB@affilsepx{\protect\\[0.5em]}
\makeatother

\author[1]{Patrick L. Combettes}
\affil[1]{North Carolina State University
\affilcr
Department of Mathematics
\affilcr
Raleigh, NC 27695, USA
\affilcr
\email{plc@math.ncsu.edu}
}
\author[2]{Diego J. Cornejo}
\affil[2]{Universidad del Pac\'{i}fico
\affilcr
Lima, Per\'u
\affilcr
\email{dj.cornejok@up.edu.pe}
}

\begin{document}
\title{%
Proximal Comixture Minimization Models for Image Recovery and 
Data Analysis\thanks{Contact author: P. L. Combettes.
Email: \email{plc@math.ncsu.edu}.
Phone: +1 919 515 2671.
This work was supported by the National
Science Foundation under grant CCF-2211123.\\
\noindent
{
A preliminary version of this work appears in the conference 
paper \cite{Eusi24} with a simpler comixture optimization model 
and without any theoretical analysis.}

}}

\date{~}

\maketitle

\vspace{12mm}

\thispagestyle{empty}
\begin{abstract}
In minimization models for image recovery and data analysis
problems, loss functions and linear operators are typically
aggregated as an average of composite terms. Each term in the
aggregate models a desired property of the ideal solution arising
from the \emph{a priori} knowledge and the observed data. We
propose an alternative minimization model based on proximal
comixtures, an operation which combines functions and linear
operators in such a way that the proximity operator of the
resulting function is computable explicitly in terms of the
individual proximity and linear operators. The mathematical
properties of this operation are analyzed and comparisons between
proximal comixtures and standard composite averages are made.
Numerical illustrations of the benefits of minimization models
based on proximal comixtures are provided in the context of image
recovery and machine learning applications.
\end{abstract}

\begin{keywords}
Convex minimization,
image recovery,
Moreau envelope,
proximal average,
proximal cocomposition,
proximal comixture.
\end{keywords}

\setcounter{page}{0}
\newpage

\section{Introduction} 

The objective of image recovery and, more broadly, of various tasks
in the areas of inverse problems and data analysis is to determine
the value of a mathematical object (an image, a signal, a set of
parameters, a distribution, a covariance matrix, a spectrum, etc.)
conveying information of interest using experimental measurements
and some \emph{a priori} knowledge (data formation model,
properties of the target solution, etc.). Over the past 60 years,
convex optimization has established itself as one of the most
reliable and efficient framework to formulate, analyze, and solve
such problems
\cite{Andr77,Cham16,Banf11,Sign21,Smms05,Cond23,Glow16,%
Judi20,Phil62,Theo20}. One typically aims at minimizing an
aggregation of convex loss functions that model individually the
desired properties of the ideal solution arising from the \emph{a
priori} knowledge and the observed data. The assumptions underlying
the minimization models to be discussed are the following (see
Section~\ref{sec:2} for notation). 

\begin{assumption}
\label{a:1}
$\HH$ is a real Hilbert space with scalar product
$\scal{\cdot}{\cdot}_{\HH}$, associated norm $\|\cdot\|_{\HH}$, and
quadratic kernel $\qq_{\HH}=\|\cdot\|_{\HH}^2/2$,
$f\in\Gamma_0(\HH)$, and $h\colon\HH\to\RR$ is a convex  and
differentiable function with a $\beta^{-1}$-Lipschitzian gradient
for some $\beta\in\RPP$. Further, $0<p\in\NN$ and, for every
$k\in\{1,\dots,p\}$, $\GG_k$ is a real Hilbert space,
$g_k\in\Gamma_0(\GG_k)$, and $0\neq L_k\colon\HH\to\GG_k$ is a
bounded linear operator. Finally, $(\alpha_k)_{1\leq k\leq p}$ are
weights in $\RPP$ which satisfy 
$\sum_{k=1}^p\alpha_k\|L_k\|^2\leq 1$. 
\end{assumption}

The most prevalent minimization framework in data analysis consists
in minimizing an objective function which aggregates the loss
functions $(g_k)_{1\leq k\leq p}$ and the linear operators
$(L_k)_{1\leq k\leq p}$ by means of a standard composite averaging
operation as follows.

\begin{problem}
\label{prob:1}
Under Assumption~\ref{a:1}, the task is to
\begin{equation}
\label{e:p1}
\minimize{x\in\HH}{f(x)+\brk1{\cav(g_k,L_k)_{1\leq k\leq p}}(x)
+h(x)},
\end{equation}
where 
\begin{equation}
\label{e:cav}
\cav(g_k,L_k)_{1\leq k\leq p}=\sum_{k=1}^p\alpha_kg_k\circ L_k
\end{equation}
is the \emph{standard composite average} of $(g_k)_{1\leq k\leq p}$
and $(L_k)_{1\leq k\leq p}$.
\end{problem}

An alternative way to aggregate the functions 
$(g_k)_{1\leq k\leq p}$ and the linear 
operators $(L_k)_{1\leq k\leq p}$ is via the
proximal comixture operation. This operation is derived from the
proximal mixture operation recently introduced in \cite{Svva23} by
duality, and it is further studied in \cite{Jota24,Eect26}. The
main objective of the present paper is to propose the use of this
aggregation process as an alternative to the standard composite
average of Problem~\ref{prob:1}. This brings us to the following
minimization model, which involves the conjugation operation of
\eqref{e:conj}. 

\begin{problem}
\label{prob:3}
Let $\gamma\in\RPP$. Under Assumption~\ref{a:1}, the task is to
\begin{equation}
\label{e:p3}
\minimize{x\in\HH}{f(x)
+\brk1{\pcmx{\gamma}(g_k,L_k)_{1\leq k\leq p}}(x)+h(x)},
\end{equation}
where
\begin{equation}
\label{e:pcm}
\pcmx{\gamma}(g_k,L_k)_{1\leq k\leq p}
=\brk4{\brk3{\sum_{k=1}^p\alpha_k\brk2{g_k^*+\gamma\qq_{\GG_k}}^*
\circ L_k}^*-\gamma\qq_{\HH}}^*
\end{equation}
is the \emph{proximal comixture} of $(g_k)_{1\leq k\leq p}$ and
$(L_k)_{1\leq k\leq p}$ with parameter $\gamma$.
\end{problem}

Let us make some observations about these two formulations to
motivate the latter.
\begin{itemize}
\item
Both Problems~\ref{prob:1} and \ref{prob:3} can be viewed as
relaxations of the ideal feasibility problem 
\begin{equation}
\label{e:feas}
\text{find}\;\;x\in\HH\quad\text{such that}\quad 
\begin{cases}
x\in\Argmin f\\
(\forall k\in\{1,\ldots,p\})\; L_kx\in\Argmin g_k\\
x\in\Argmin h,
\end{cases}
\end{equation}
which aims at imposing all the desired properties exactly and hence
at minimizing all the loss functions simultaneously. However, the
set 
\begin{equation}
\label{e:feasZ}
Z=\Argmin f\cap\brk3{\bigcap_{k=1}^pL_k^{-1}(\Argmin g_k)}
\cap\Argmin h
\end{equation}
of solutions to \eqref{e:feas} is typically empty and, as will be
seen in Remark~\ref{r:35}\ref{r:35ii}, both Problems~\ref{prob:1}
and \ref{prob:3} are relaxations of \eqref{e:feas} in the sense
that, if $Z\neq\emp$, then $Z$ is the set of solutions of these 
two problems. 
\item
The aggregation in Problem~\ref{prob:1} may not be robust to
perturbations. For instance, let us consider the special case when
$f=0$, $h=0$, and, for every $k\in\{1,\ldots,p\}$, $\GG_k=\HH$,
$L_k=\Id_{\HH}$, and $g_k$ is the 
indicator function of a nonempty closed convex set
$C_k\subset\HH$. This reduces \eqref{e:p1} to 
\begin{equation}
\label{e:p12}
\text{find}\;\;x\in\bigcap_{k=1}^p C_k,
\end{equation}
which happens to coincide with \eqref{e:feas} in this case.
If the sets $(C_k)_{1\leq k\leq p}$ are not specified exactly, no
solution may exist \cite{Cens18,Sign94}. However, it will follow
from Remark~\ref{r:35}\ref{r:35i} that the solutions to
Problem~\ref{prob:3} minimize the average squared distance to the
sets. Such solutions are classical surrogates in inconsistent
feasibility problems \cite{Sign94}, which go back to Legendre's
least-squares method for systems of linear equations \cite{Lege05}.
They exist under mild conditions, such as the boundedness of one
of the sets.
\item
Problem~\ref{prob:1} involves a simple aggregation process by
averaging. However, the nonsmooth function \eqref{e:cav} has no
explicit proximity operator, and solving \eqref{e:p1} therefore
calls for sophisticated proximal splitting methods to decompose
the $p$ functions $(g_k)_{1\leq k\leq p}$ and the $p$ linear
operators $(L_k)_{1\leq k\leq p}$ individually \cite{Acnu24}.
Overall, solving Problem~\ref{prob:1} requires the splitting of
$2p+2$ terms, which can be expected to lead to algorithms that are
slower and necessitate more memory storage than those that would
split less terms. By contrast, the aggregated function
\eqref{e:pcm} in Problem~\ref{prob:3} is less intuitive but, as
will be seen in Proposition~\ref{p:1}\ref{p:1ii}, its proximity
operator can be computed explicitly in terms of those of the
functions $(g_k)_{1\leq k\leq p}$ and of the linear operators
$(L_k)_{1\leq k\leq p}$ as
\begin{equation}
\label{e:jjm1}
\prox_{\gamma\pcmx{\gamma}(g_k,L_k)_{1\leq k\leq p}}
=\Id_{\HH}-\sum_{k=1}^p\alpha_kL_k^*\circ\brk1{\Id_{\GG_k}
-\prox_{\gamma g_k}}\circ L_k.
\end{equation}
In turn, solving Problem~\ref{prob:3} requires the splitting of
only 3 terms.
\item
In the special case when $\sum_{k=1}^p\alpha_k=1$ and, for every
$k\in\{1,\ldots,p\}$, $\GG_k=\HH$ and $L_k=\Id_{\HH}$, the proximal
comixture of \eqref{e:pcm} reduces to the proximal average 
\begin{equation}
\label{e:pav}
\pave{\gamma}(g_k)_{1\leq k\leq p}
=\brk4{\brk3{\sum_{k=1}^p\alpha_k\brk2{
g_k^*+\gamma\qq_{\HH}}^*}^*-\gamma\qq_{\HH}}^*,
\end{equation}
and \eqref{e:jjm1} gives
\begin{equation}
\label{e:jjm2}
\prox_{\gamma\pave{\gamma}(g_k)_{1\leq k\leq p}}
=\sum_{k=1}^p\alpha_k\prox_{\gamma g_k}.
\end{equation}
This aggregation process has been implicitly introduced in
\cite{More65}, extensively studied in \cite{Baus08}, and used in
data analysis problems in 
\cite{Ang25,Cheu17,Huzz22,Liti20,Shen17,Yuyl13}, where its benefits
over the standard average
\begin{equation}
\label{e:avg}
\ave(g_k)_{1\leq k\leq p}=\sum_{k=1}^p\alpha_kg_k,
\end{equation}
are discussed. 
\item
Formally, evaluating \eqref{e:pcm} at $\gamma=0$ gives back
\eqref{e:cav}. This connection will be made mathematically precise
in Theorem~\ref{t:4}, where we study the asymptotic behavior as
$\gamma\downarrow 0$.
\end{itemize}

Notation and background are provided in Section~\ref{sec:2}. In
Section~\ref{sec:3}, we study the mathematical properties of
proximal comixtures and investigate connections between
Problems~\ref{prob:1} and \ref{prob:3}. Section~\ref{sec:4} is
devoted to algorithms for solving Problems~\ref{prob:1} and
\ref{prob:3}. Finally, Section~\ref{sec:5} provides numerical
illustrations of proximal comixture models in concrete signal
restoration, image reconstruction, and linear regression problems. 

\section{Notation and background}
\label{sec:2}
Throughout, $\HH$ is a real Hilbert space with power set $2^{\HH}$,
identity operator $\Id_{\HH}$, scalar product
$\scal{\Cdot}{\Cdot}_{\HH}$, associated norm $\norm{\Cdot}_{\HH}$,
and quadratic kernel $\qq_{\HH}=\norm{\Cdot}_{\HH}^2/2$. For
background on nonlinear analysis in Hilbert spaces, see
\cite{Livre1}. 

Let $\varphi\colon\HH\to\RXX$. The conjugate of $\varphi$ is
\begin{equation}
\label{e:conj}
\varphi^*\colon\HH\to\RXX\colon
x^*\mapsto\sup_{x\in\HH}\brk1{\scal{x}{x^*}_{\HH}-\varphi(x)}.
\end{equation}
If $\inf_{y\in\HH}\varphi(y)>\minf$ and $\varepsilon\in\RP$, the
set of $\varepsilon$-minimizers of $\varphi$ is 
\begin{equation}
\varepsilon\text{--}\!\Argmin\varphi=\menge{x\in\HH}{
\varphi(x)\leq\inf_{y\in\HH}\varphi(y)+\varepsilon}.
\end{equation}
In particular, the set of minimizers of $\varphi$ is
$\Argmin\varphi=0\text{--}\!\Argmin\varphi$.
If $\Argmin\varphi$ is a singleton, its unique element is denoted
by $\argmin_{x\in\HH}\varphi(x)$.
Moreover, $\varphi$ is proper if $\minf\notin\varphi(\HH)$ and
$\dom\varphi=\menge{x\in\HH}{\varphi(x)<\pinf}\neq\emp$. If
$\varphi$ is proper, its subdifferential is
\begin{equation}
\partial \varphi\colon\HH\to 2^\HH\colon x\mapsto
\menge{x^*\in\HH}{(\forall y\in\HH)\;
\scal{y-x}{x^*}_{\HH}+\varphi(x)\leq \varphi(y)}.
\end{equation}
Let $\gamma\in\RPP$. Then the \emph{(lower) Moreau envelope} of
$\varphi$ with parameter $\gamma$ is 
\begin{equation}
\label{e:m1}
\lenv{\varphi}{\gamma}\colon\HH\to\RXX\colon x\mapsto
\inf_{y\in\HH}\brk2{\varphi(y)+\dfrac{1}{\gamma}\qq_{\HH}(x-y)},
\end{equation}
and the \emph{upper Moreau envelope} of $\varphi$ with parameter
$\gamma$ is 
\begin{equation}
\label{e:m2}
\uenv{\varphi}{\gamma}\colon\HH\to\RXX\colon x\mapsto
\sup_{y\in\HH}\brk2{\varphi(y)-\dfrac{1}{\gamma}\qq_{\HH}(x-y)}.
\end{equation}
Given $\rho\in\RR$, $\Gamma_{\rho}(\HH)$ denotes the class of
proper lower semicontinuous functions $\varphi\colon\HH\to\RX$
such that $\varphi+\rho\qq_{\HH}$ is convex. The proximity operator
of $\varphi\in\Gamma_0(\HH)$ is
\begin{equation}
\label{e:prox}
\prox_{\varphi}\colon\HH\to\HH\colon x\mapsto
\argmin_{y\in\HH}\brk1{\varphi(y)+\qq_{\HH}(x-y)}
\end{equation}
and it is characterized by
\begin{equation}
\label{e:jjm}
(\forall x\in\HH)(\forall p\in\HH)\quad
p=\prox_{\varphi}x\quad\Leftrightarrow\quad x-p\in
\partial\varphi(p).
\end{equation}
An operator $T\colon\HH\to\HH$ is $\gamma$-cocoercive if
\begin{equation}
\label{e:coco}
(\forall x\in\HH)(\forall y\in\HH)\quad
\scal{x-y}{Tx-Ty}_{\HH}\geq\gamma\|Tx-Ty\|_{\HH}^2.
\end{equation}
Let $C\subset\HH$. Then the indicator function of $C$ is denoted by
$\iota_C$ and the distance function to $C$ is denoted by $d_C$. If
$C$ is nonempty, closed, and convex, its projection operator is
denoted by $\proj_C$. The closed ball with center $x\in\HH$ and
radius $\rho\in\RPP$ is denoted by $B(x;\rho)$.

The following results will be useful.

\begin{lemma}
\label{l:2}
Let $\varphi\in\Gamma_0(\HH)$ and let $\gamma\in\RPP$. Then the
following hold:
\begin{enumerate}
\item
\label{l:2i}
$\lenv{\varphi}{\gamma}=\brk1{\varphi^*+\gamma\qq_{\HH}}^*$.
\item
\label{l:2ii}
$\uenv{\varphi}{\gamma}=\brk1{\varphi^*-\gamma\qq_{\HH}}^*$.
\end{enumerate}
\end{lemma}
\begin{proof}
Recall that, since $\varphi\in\Gamma_0(\HH)$, $\varphi^{**}
=\varphi$ \cite[Corollary~13.38]{Livre1}.

\ref{l:2i}: This follows from \cite[Proposition~14.1]{Livre1}.

\ref{l:2ii}: Apply \cite[Theorem~2.2]{Hiri86} with
$g=\varphi^*$ and $h=\gamma\qq_{\HH}$. 
\end{proof}

We recall the following result from \cite[Lemma~3]{Bern10} 
(see also \cite[Example~11.26(d)]{Rock09} for the
finite-dimensional case).

\begin{lemma}
\label{l:4}
Let $\varphi\colon\HH\to\RX$ be proper and let $\gamma\in\RPP$.
Then the following are equivalent:
\begin{enumerate}
\item
\label{l:4i}
$\uenv{\brk{\lenv{\varphi}{\gamma}}}{\gamma}=\varphi$.
\item
\label{l:4ii}
$\varphi\in\Gamma_{1/\gamma}(\HH)$.
\end{enumerate}
\end{lemma}

\begin{lemma}
\label{l:6}
Let $\varphi\colon\HH\to\RR$ be continuous and convex, and let
$\gamma\in\RPP$. Then the following are equivalent:
\begin{enumerate}
\item
\label{l:6i}
$\lenv{\brk{\uenv{\varphi}{\gamma}}}{\gamma}=\varphi$. 
\item
\label{l:6ii}
$-\varphi\in\Gamma_{1/\gamma}(\HH)$.
\item
\label{l:6iii}
$\brk{\uenv{\varphi}{\gamma}}^*=\varphi^*-\gamma\qq_{\HH}$.
\item
\label{l:6iv}
$\uenv{\varphi}{\gamma}\in\Gamma_0(\HH)$,
$\varphi$ is Fr\'{e}chet differentiable, and
$\prox_{\gamma\brk{\uenv{\varphi}{\gamma}}}
=\Id_{\HH}-\gamma\nabla \varphi$.
\item
\label{l:6v}
$\varphi$ is Fr\'{e}chet differentiable and
$\nabla\varphi$ is $\gamma$-cocoercive.
\item
\label{l:6vi}
$\varphi$ is Fr\'{e}chet differentiable and
$\nabla\varphi$ is $\gamma^{-1}$-Lipschitzian.
\end{enumerate}
\end{lemma}
\begin{proof}
It follows from \eqref{e:m1} and \eqref{e:m2} that
$\lenv{\varphi}{\gamma}=-\uenv{(-\varphi)}{\gamma}$. Therefore,
\begin{equation}
\label{e:l6}
\lenv{\brk1{\uenv{\varphi}{\gamma}}}{\gamma}
=-\uenv{\brk1{-\uenv{\varphi}{\gamma}}}{\gamma}
=-\uenv{\brk1{\lenv{(-\varphi)}{\gamma}}}{\gamma}.
\end{equation}

\ref{l:6i}$\Rightarrow$\ref{l:6ii}: 
By \eqref{e:l6}, 
$\uenv{\brk{\lenv{(-\varphi)}{\gamma}}}{\gamma}=-\varphi$. Thus,
Lemma~\ref{l:4} yields $-\varphi\in\Gamma_{1/\gamma}(\HH)$. 

\ref{l:6ii}$\Rightarrow$\ref{l:6iii}:
By \cite[Proposition~14.2]{Livre1},
$\varphi^*\in\Gamma_{-\gamma}(\HH)$. Thus, the
result follows from Lemma~\ref{l:2}\ref{l:2ii} and
\cite[Corollary~13.38]{Livre1}.

\ref{l:6iii}$\Rightarrow$\ref{l:6iv}: 
Since $\varphi\in\Gamma_0(\HH)$, \cite[Corollary~13.38]{Livre1}
asserts that $\varphi^*\in\Gamma_0(\HH)$, which implies that
$\uenv{\varphi}{\gamma}$ is proper. In turn, it results from 
\cite[Proposition~13.13]{Livre1} and Lemma~\ref{l:2}\ref{l:2ii}
that $\uenv{\varphi}{\gamma}\in\Gamma_0(\HH)$. Thus, we derive
from Lemma~\ref{l:2}\ref{l:2i} and \cite[Corollary~13.38]{Livre1}
that 
\begin{equation}
\lenv{\brk1{\uenv{\varphi}{\gamma}}}{\gamma}
=\brk2{\brk1{\uenv{\varphi}{\gamma}}^*+\gamma\qq_{\HH}}^*
=\brk2{\brk1{\varphi^*-\gamma\qq_{\HH}}+\gamma\qq_{\HH}}^*
=\varphi^{**}=\varphi.
\end{equation}
Hence, \cite[Proposition~12.30]{Livre1} guarantees that 
$\varphi=\lenv{\brk1{\uenv{\varphi}{\gamma}}}{\gamma}$ is
Fr\'{e}chet differentiable and that
$\nabla\varphi=\gamma^{-1}\brk{\Id_{\HH}
-\prox_{\gamma\brk{\uenv{\varphi}{\gamma}}}}$.

\ref{l:6iv}$\Rightarrow$\ref{l:6v}:
This follows from \cite[Proposition~12.28]{Livre1}.

\ref{l:6v}$\Rightarrow$\ref{l:6vi}:
This follows from \eqref{e:coco} and the Cauchy--Schwarz
inequality.

\ref{l:6vi}$\Rightarrow$\ref{l:6i}:
It follows from the equivalence (i)$\Leftrightarrow$(vi) in
\cite[Theorem~18.15]{Livre1} that
$-\varphi\in\Gamma_{1/\gamma}(\HH)$.
We therefore deduce from Lemma~\ref{l:4} that
$\uenv{\brk{\lenv{(-\varphi)}{\gamma}}}{\gamma}=-\varphi$.
We conclude via \eqref{e:l6}.
\end{proof}

\begin{remark}
The functions $\lenv{\varphi}{\gamma}$ and 
$\uenv{\varphi}{\gamma}$ are, respectively, the infimal and the
supremal convolutions of $\varphi$ and $\qq_{\HH}/\gamma$.
In \cite[Definition~2.4]{Hiri94}, 
$\uenv{\varphi}{\gamma}$ is called the deconvolution of 
$\varphi$ by $\qq_{\HH}/\gamma$. Furthermore,
when $\gamma<\rho$, the functions
$\lenv{\brk{\uenv{\varphi}{\rho}}}{\gamma}$ and
$\uenv{\brk{\lenv{\varphi}{\rho}}}{\gamma}$ are known as the
\emph{Lasry-Lions regularizations} of $\varphi$, which were
introduced in \cite{Lasry86} and further studied in
\cite{Atto93,Bern10,Rock09,Strom96}, whereas 
$\lenv{\brk{\uenv{\varphi}{\gamma}}}{\gamma}$ is called
\emph{proximal hull} of $\varphi$ in \cite[Example~1.44]{Rock09}.
\end{remark}

\section{Properties of proximal comixtures}
\label{sec:3}

Throughout this section, Assumption~\ref{a:1} is in force. Recall
from \eqref{e:cav} that the associated standard composite average
is
\begin{equation}
\label{e:cav2}
\cav(g_k,L_k)_{1\leq k\leq p}=\sum_{k=1}^p\alpha_kg_k\circ L_k,
\end{equation}
and from \eqref{e:pcm} that the associated proximal comixture with
parameter $\gamma\in\RPP$ is
\begin{equation}
\label{e:pcm2}
\pcmx{\gamma}(g_k,L_k)_{1\leq k\leq p}
=\brk4{\brk3{\sum_{k=1}^p\alpha_k\brk2{
g_k^*+\gamma\qq_{\GG_k}}^*\circ L_k}^*-\gamma\qq_{\HH}}^*.
\end{equation}

\begin{figure}[t]
\centering
\vskip 20mm
\scalebox{0.6} 
{
\begin{pspicture}(10.0,-1.7)(30.0,2.0)

\psset{unit=1cm} 
\rput[l](7,0){ 
\rput(0.1,0.0){\large$\boldsymbol{{g_1}}$}
\psline[linewidth=0.04cm,arrowsize=2.2mm]{->}(0.5,0.0)(1.5,0.0)
\psframe[linewidth=0.04,dimen=outer,linecolor=dblue,framearc=0.4]%
(1.5,0.5)(2.5,-0.5)
\psline[linewidth=0.04cm,arrowsize=2.2mm,linecolor=dblue]{->}%
(2.5,0.0)(3.5,0.0)
\psline[linewidth=0.04cm,arrowsize=2.2mm,linecolor=dblue]{->}%
(2.0,1.5)(2.0,0.5)
\rput(0.1,-2){\large$\boldsymbol{{g_p}}$}
\psline[linewidth=0.04cm,arrowsize=2.2mm]{->}(0.5,-2)(1.5,-2)
\psframe[linewidth=0.04,dimen=outer,linecolor=dblue,framearc=0.4]%
(1.5,-1.5)(2.5,-2.5)
\psline[linewidth=0.04cm,arrowsize=2.2mm,linecolor=dblue]{->}%
(2.5,-2)(3.5,-2)
\psline[linewidth=0.04cm,arrowsize=2.2mm,linecolor=dblue]{->}%
(2.0,-3.5)(2.0,-2.5)

\psframe[linewidth=0.04,dimen=outer,linecolor=dblue,framearc=0.4]%
(3.5,0.45)(5.5,-2.45)
\psline[linewidth=0.04cm,arrowsize=2.2mm]{->}(5.5,-1)(6.5,-1)
\rput(0.1,-0.9){\huge$\boldsymbol{{\vdots}}$}

\psframe[linewidth=0.04,dimen=outer,linecolor=dblue,framearc=0.4,%
linestyle=dashed](1.0,2.4)(6.1,-4.4)

\rput(8.3,-1){\Large${\cav(g_k,L_k)_{1\leq k\leq p}}$}

\textcolor{dblue}{
\rput(2.9,-0.9){\huge$\boldsymbol{{\vdots}}$}
\rput(2.0,0.0){\Large$\boldsymbol{\circ}$}
\rput(2.0,-2){\Large$\boldsymbol{\circ}$}
\rput(2.9,0.2){$\boldsymbol{\alpha_1}$}
\rput(2.9,-2.3){$\boldsymbol{\alpha_p}$}
\rput(2.0,1.8){\large$\boldsymbol{L_1}$}
\rput(2.0,-3.8){\large$\boldsymbol{L_p}$}
\rput(4.5,-1){\huge$\boldsymbol{\sum}$}
\rput(3.5,3.5){\textrm{\parbox[c]{2cm}{\centering
\textbf{standard}\\\textbf{composite}\\\textbf{average}}}}
}

}
\rput[r](19.,0){ 
\rput(0.1,0.0){\large$\boldsymbol{g_1}$}
\psline[linewidth=0.04cm,arrowsize=2.2mm,linecolor=black]{->}%
(0.5,0.0)(1.5,0.0)
\psframe[linewidth=0.04,dimen=outer,linecolor=dgreen,framearc=0.4]%
(1.5,0.6)(3.5,-0.6)
\psline[linewidth=0.04cm,arrowsize=2.2mm,linecolor=dgreen]{->}%
(3.5,0.0)(4.5,0.0)

\rput(0.1,-2){\large$\boldsymbol{g_p}$}
\psline[linewidth=0.04cm,arrowsize=2.2mm,linecolor=black]{->}%
(0.5,-2)(1.5,-2)
\psframe[linewidth=0.04,dimen=outer,linecolor=dgreen,framearc=0.4]%
(1.5,-1.4)(3.5,-2.6)
\psline[linewidth=0.04cm,arrowsize=2.2mm,linecolor=dgreen]{->}%
(3.5,-2.0)(4.5,-2.0)

\psline[linewidth=0.04cm,arrowsize=2.2mm,linecolor=dred]{->}%
(6.5,-1)(7.5,-1)
\rput(0.1,-0.9){\huge$\boldsymbol{{\vdots}}$}

\psframe[linewidth=0.04,dimen=outer,linecolor=dred,framearc=0.4]%
(7.5,-0.4)(9.5,-1.6)
\psframe[linewidth=0.04,dimen=outer,linecolor=dred,framearc=0.4,%
linestyle=dashed](0.9,1.4)(10.1,-3.4)
\psline[linewidth=0.04cm,arrowsize=2.2mm]{->}(9.5,-1)(10.5,-1)

\rput(12.5,-1.05){\Large${\pcmx{\gamma}(g_k,L_k)_{1\leq k\leq p}}$}

\psframe[linewidth=0.04,dimen=outer,linecolor=dblue,framearc=0.4,
linestyle=dashed](4.5,0.55)(6.5,-2.55)
\textcolor{dblue}{
\rput(5.5,-1){\textrm{\parbox[c]{2cm}{
\centering\textbf{standard}\\
\textbf{composite}\\
\textbf{average}}}}}

\textcolor{dgreen}{
\rput(3.8,-0.9){\huge$\boldsymbol{{\vdots}}$}
\rput(2.4,0.0){\large$\boldsymbol{\lenv{}{\gamma}}$}
\rput(2.4,-2.0){\large$\boldsymbol{\lenv{}{\gamma}}$}
}

\textcolor{dred}{
\rput(8.35,-1){\large$\boldsymbol{\uenv{}{\gamma}}$}
\rput(5.5,2.4){\textrm{\parbox[c]{2cm}%
{\centering\textbf{proximal}\\
\textbf{comixture}}}}}
}

\end{pspicture}
}
\vskip 20mm
\caption{
(\textcolor{dblue}{left}): Standard composite average
\eqref{e:cav2}.
(\textcolor{dred}{right}): Proximal comixture \eqref{e:pcm2} in
terms of the Moreau envelopes of \eqref{e:m1} and \eqref{e:m2}
using Proposition~\ref{p:0}\ref{p:0i}.
}
\label{fig:1}
\end{figure}

We first provide reformulations of the proximal comixture of
\eqref{e:pcm2} in terms of the Moreau envelopes of \eqref{e:m1} and
\eqref{e:m2} (see Fig.~\ref{fig:1}), and then in terms of the 
proximal average of \eqref{e:pav} and of the following proximal 
cocomposition operation.
These connections will not only shed new light on proximal
comixtures but also play a role in forthcoming proofs.

\begin{definition}[\protect{\cite{Svva23}}]
\label{d:pcc}
Let $\GG$ be a real Hilbert space, let $g\in\Gamma_0(\GG)$, and let
$L\colon\HH\to\GG$ be a bounded linear operator. The 
\emph{proximal cocomposition} of $g$ and $L$ with parameter
$\gamma\in\RPP$ is
\begin{equation}
\label{e:r1b}
L\proxcc{\gamma}g=\brk3{\brk2{
\brk1{g^*+\gamma\qq_{\GG}}^*\circ L}^*-\gamma\qq_{\HH}}^*.
\end{equation}
\end{definition}

\begin{proposition}
\label{p:0}
Let $\gamma\in\RPP$. Then the following hold: 
\begin{enumerate}
\item
\label{p:0i}
$\pcmx{\gamma}(g_k,L_k)_{1\leq k\leq p}=\uenv{
\sum_{k=1}^p\alpha_k\brk{\lenv{g_k}{\gamma}}\circ L_k}{\gamma}$.
\item
\label{p:0ii}
Suppose that $\sum_{k=1}^p\alpha_k=1$ and 
$\max_{1\leq k\leq p}\|L_k\|\leq 1$. Then
$\pcmx{\gamma}(g_k,L_k)_{1\leq k\leq p}=
\pave{\gamma}\brk{L_k\proxcc{\gamma}g_k}_{1\leq k\leq p}$.
\item
\label{p:0iii}
Let $\GG$ be the standard product vector space 
$\bigtimes_{k=1}^p\GG_k$, with generic element
$\boldsymbol{y}=(y_k)_{1\leq k\leq p}$, and equipped with the scalar
product 
\begin{equation}
\label{e:sp}
\scal{\cdot}{\cdot}_{\GG}\colon(\boldsymbol{y},\boldsymbol{y}')
\mapsto\sum_{k=1}^p\alpha_k\scal{y_k}{y_k'}_{\GG_k}, 
\end{equation}
and set
\begin{equation}
\label{e:p0iii}
L\colon\HH\to\GG\colon x\mapsto(L_kx)_{1\leq k\leq p}
\quad\text{and}\quad
g\colon\GG\to\RX\colon\boldsymbol{y}\mapsto
\sum_{k=1}^p\alpha_kg_k(y_k).
\end{equation}
Then $\pcmx{\gamma}(g_k,L_k)_{1\leq k\leq p}=L\proxcc{\gamma}g$.
\end{enumerate}
\end{proposition}
\begin{proof}
Set 
$\varphi=\sum_{k=1}^p\alpha_k\brk{\lenv{g_k}{\gamma}}\circ L_k$.
By virtue of \cite[Proposition~12.30]{Livre1},
\begin{equation}
\label{e:0}
\nabla\varphi=\dfrac{1}{\gamma}\sum_{k=1}^p
\alpha_kL_k^*\circ(\Id_{\GG_k}-\prox_{\gamma g_k})\circ L_k.
\end{equation} 
However, $\sum_{k=1}^p\alpha_k\|L_k\|^2\leq 1$ by
Assumption~\ref{a:1} and \cite[Proposition~12.28]{Livre1} asserts 
that the operators 
$(\Id_{\GG_k}-\prox_{\gamma g_k})_{1\leq k\leq p}$ are 
nonexpansive. Therefore, \eqref{e:0} implies that $\nabla\varphi$
is $\gamma^{-1}$-Lipschitzian.

\ref{p:0i}:
We derive from \eqref{e:pcm2} and 
Lemma~\ref{l:2}\ref{l:2i}--\ref{l:2ii} that 
\begin{equation}
\pcmx{\gamma}(g_k,L_k)_{1\leq k\leq p}
=\brk1{\varphi^*-\gamma\qq_{\HH}}^*
=\uenv{\varphi}{\gamma}.
\end{equation}

\ref{p:0ii}:
It follows from \cite[Proposition~3.13(ii)]{Eect26} that
$(\forall k\in\{1,\ldots,p\})$ 
$\lenv{\brk{L_k\proxcc{\gamma}g_k}}{\gamma}
=\brk{\lenv{g_k}{\gamma}}\circ L_k$.
Therefore, \ref{p:0i}, Lemma~\ref{l:2}\ref{l:2i}--\ref{l:2ii}, and
\eqref{e:pav} yield
\begin{align}
\pcmx{\gamma}(g_k,L_k)_{1\leq k\leq p}
&=\uenv{\varphi}{\gamma}\nonumber\\
&=\uenv{\brk4{\sum_{k=1}^p\alpha_k
\lenv{\brk1{L_k\proxcc{\gamma}g_k}}{\gamma}}}{\gamma}\nonumber\\
&=\brk4{\brk3{\displaystyle\sum_{k=1}^p\alpha_k
\lenv{\brk1{L_k\proxcc{\gamma}g_k}}{\gamma}}^*
-\gamma\qq_{\HH}}^*\nonumber\\
&=\pave{\gamma}\brk1{L_k\proxcc{\gamma}g_k}_{1\leq k\leq p}.
\end{align}

\ref{p:0iii}: Since $\|\cdot\|_{\GG}^2\colon\GG\to\RR\colon
\boldsymbol{y}\mapsto\sum_{k=1}^p\alpha_k\|y_k\|_{\GG_k}^2$,
$\lenv{g}{\gamma}\colon\GG\to\RR\colon\boldsymbol{y}
\mapsto\sum_{k=1}^p\alpha_k(\lenv{g_k}{\gamma})(y_k)$.
Therefore, $\brk{\lenv{g}{\gamma}}\circ L
=\sum_{k=1}^p\alpha_k\brk{\lenv{g_k}{\gamma}}\circ L_k$, 
and the result follows from \ref{p:0i}, items \ref{l:2ii} and
\ref{l:2i} in Lemma~\ref{l:2}, and \eqref{e:r1b}.
\end{proof}

\begin{proposition}
\label{p:1}
Let $\gamma\in\RPP$. Then the following hold: 
\begin{enumerate}
\item
\label{p:1i}
$\pcmx{\gamma}(g_k,L_k)_{1\leq k\leq p}\in\Gamma_0(\HH)$.
\item
\label{p:1ii}
$\prox_{\gamma\pcmx{\gamma}(g_k,L_k)_{1\leq k\leq p}}
=\Id_{\HH}-\sum_{k=1}^p\alpha_kL_k^*\circ
(\Id_{\GG_k}-\prox_{\gamma g_k})\circ L_k$.
\item
\label{p:1iii}
Suppose that one of the following is satisfied:
\begin{enumerate}
\item
\label{p:1iiia}
$\sum_{k=1}^p\alpha_k\|L_k\|^2<1$.
\item
\label{p:1iiia+}
For every $k\in\{1,\ldots,p\}$, $\dom g_k=\GG_k$.
\item
\label{p:1iiib}
$\sum_{k=1}^p\alpha_k=1$, $\max_{1\leq k\leq p}\|L_k\|\leq 1$, and
there exists ${\ell}\in\{1,\ldots,p\}$ such that $\dom
g_{\ell}=\GG_{\ell}$.
\end{enumerate}
Then $\dom\pcmx{\gamma}(g_k,L_k)_{1\leq k\leq p}=\HH$.
\item
\label{p:1iv}
$\lenv{\pcmx{\gamma}(g_k,L_k)_{1\leq k\leq p}}{\gamma}
=\sum_{k=1}^p\alpha_k\brk{\lenv{g_k}{\gamma}}\circ L_k$.
\item
\label{p:1v}
$\Argmin\pcmx{\gamma}(g_k,L_k)_{1\leq k\leq p}
=\Argmin\sum_{k=1}^p\alpha_k\brk{\lenv{g_k}{\gamma}}\circ L_k$.
\item
\label{p:1vi}
Suppose that $\bigcap_{k=1}^pL_k^{-1}(\Argmin g_k)\neq\emp$. Then 
\begin{equation}
\label{e:24}
\Argmin\pcmx{\gamma}(g_k,L_k)_{1\leq k\leq p}
=\Argmin\cav(g_k,L_k)_{1\leq k\leq p}
=\bigcap_{k=1}^pL_k^{-1}(\Argmin g_k).
\end{equation}
\end{enumerate}
\end{proposition}
\begin{proof}
Define $\GG$, $L$, and $g$ as in \eqref{e:p0iii}, and note that
\begin{equation}
\label{e:prop1}
(\forall x\in\HH)\quad
\|Lx\|_{\GG}^2
=\sum_{k=1}^p\alpha_k\|L_kx\|^2_{\GG_k}
\leq\sum_{k=1}^p\alpha_k\|L_k\|^2\|x\|^2_{\HH},
\end{equation}
which yields $\|L\|^2\leq\sum_{k=1}^p\alpha_k\|L_k\|^2\leq 1$ by
Assumption~\ref{a:1}. Further, by 
Proposition~\ref{p:0}\ref{p:0iii},
\begin{equation}
\label{e:lien}
\pcmx{\gamma}(g_k,L_k)_{1\leq k\leq p}=L\proxcc{\gamma}g. 
\end{equation}

\ref{p:1i}: This follows from \eqref{e:lien} and 
\cite[Proposition~3.7(i)]{Eect26}.

\ref{p:1ii}: Note that
$(\forall (y_k^*)_{1\leq k\leq p}\in\GG)$ 
$L^*(y_k^*)_{1\leq k\leq p}=\sum_{k=1}^p\alpha_kL_k^*y_k^*$ and 
$\prox_{\gamma g}(y_k^*)_{1\leq k\leq p}
=(\prox_{\gamma g_k}y_k^*)_{1\leq k\leq p}$.
It therefore follows from \eqref{e:lien} and
\cite[Proposition~3.10(ii)]{Eect26} that
\begin{align}
\prox_{\gamma\pcmx{\gamma}(g_k,L_k)_{1\leq k\leq p}}
&=\prox_{\gamma(L\proxcc{\gamma}g)}\nonumber\\
&=\Id_{\HH}-L^*\circ(\Id_{\GG}-\prox_{\gamma g})\circ L
\nonumber\\
&=\Id_{\HH}-\sum_{k=1}^p\alpha_kL_k^*\circ
(\Id_{\GG_k}-\prox_{\gamma g_k})\circ L_k.
\end{align}

\ref{p:1iiia}:
In this case, $\|L\|^2\leq\sum_{k=1}^p\alpha_k\|L_k\|^2<1$.
Therefore, the assertion follows from \eqref{e:lien} and
\cite[Proposition~3.2(iv)(a)]{Eect26}.

\ref{p:1iiia+}:
Since $\dom g=\bigtimes_{k=1}^p\dom g_k
=\bigtimes_{k=1}^p\GG_k=\GG$,
the assertion follows from \eqref{e:lien} and
\cite[Proposition~3.2(iv)(b)]{Eect26}.

\ref{p:1iiib}:
According to \cite[Proposition~3.2(iv)(b)]{Eect26}, 
$\dom\brk{L_{\ell}\proxcc{\gamma}g_{\ell}}=\HH$. Hence, we derive
from Proposition~\ref{p:0}\ref{p:0ii} and
\cite[Remark~5.11(v)]{Svva23} that
\begin{equation}
\dom\pcmx{\gamma}(g_k,L_k)_{1\leq k\leq p}
=\dom\pave{\gamma}\brk1{L_k\proxcc{\gamma}g_k}_{1\leq k\leq p}
=\sum_{k=1}^p\alpha_k\dom\brk1{L_k\proxcc{\gamma}g_k}
=\HH.
\end{equation}

\ref{p:1iv}:
Set 
$\varphi=\sum_{k=1}^p\alpha_k\brk{\lenv{g_k}{\gamma}}\circ L_k$.
As seen after \eqref{e:0}, $\nabla\varphi$ is
$\gamma^{-1}$-Lipschitzian. On the other hand, 
Proposition~\ref{p:0}\ref{p:0i} asserts that 
$\pcmx{\gamma}(g_k,L_k)_{1\leq k\leq p}=\uenv{\varphi}{\gamma}$.
Altogether, appealing to the equivalence 
\ref{l:6i}$\Leftrightarrow$\ref{l:6vi} in Lemma~\ref{l:6}, we
conclude that 
$\lenv{\pcmx{\gamma}(g_k,L_k)_{1\leq k\leq p}}{\gamma}
=\lenv{(\uenv{\varphi}{\gamma})}{\gamma}=\varphi$.

\ref{p:1v}:
The result follows from \ref{p:1i}, \ref{p:1iv}, and the fact that
the set of minimizers of a function in $\Gamma_0(\HH)$ coincides
with that of its lower Moreau envelope
\cite[Proposition~17.5]{Livre1}.

\ref{p:1vi}:
Since $\bigcap_{k=1}^pL_k^{-1}(\Argmin g_k)\neq\emp$,
\cite[Proposition~17.5]{Livre1} and \ref{p:1v} imply that
\begin{align}
\label{e:p35b}
\Argmin\cav(g_k,L_k)_{1\leq k\leq p}
&=\bigcap_{k=1}^pL_k^{-1}(\Argmin g_k)\nonumber\\
&=\bigcap_{k=1}^pL_k^{-1}(\Argmin\lenv{g_k}{\gamma})\nonumber\\
&=\Argmin\sum_{k=1}^p\alpha_k\brk{\lenv{g_k}{\gamma}}\circ L_k
\nonumber\\
&=\Argmin\pcmx{\gamma}(g_k,L_k)_{1\leq k\leq p},
\end{align}
which completes the proof.
\end{proof}

\begin{remark}
\label{r:35}
Let us make a couple of observations about 
Proposition~\ref{p:1}.
\begin{enumerate}
\item
\label{r:35i}
Suppose that, in Assumption~\ref{a:1}, $f=h=0$ and, for every
$k\in\{1,\ldots,p\}$, $g_k=\iota_{D_k}$ for some nonempty closed
convex set $D_k\subset\GG_k$. Then Problem~\ref{prob:1} amounts to
finding $x\in\HH$ such that, for every $k\in\{1,\ldots,p\}$, 
$L_kx\in D_k$. On the other hand, it follows from 
Proposition~\ref{p:1}\ref{p:1v} and \cite[Example~12.21]{Livre1}
that Problem~\ref{prob:3} amounts to finding a minimizer of the
least-squares function 
$x\mapsto\sum_{k=1}^p\alpha_kd_{D_k}^2(L_kx)$, which has been used
in the relaxation of inconsistent feasibility problems
\cite{Cens05,Sign94,Siim19}.
\item
\label{r:35ii}
In general, the solution sets of Problems~\ref{prob:1} and
\ref{prob:3} are distinct. However, it follows from 
Proposition~\ref{p:1}\ref{p:1vi} that, if the set $Z$ in
\eqref{e:feasZ} is nonempty, then it coincides with the solution
set of both problems.
\end{enumerate}
\end{remark}

Next, we compare the standard composite average of
$\cav(g_k,L_k)_{1\leq k\leq p}$ of \eqref{e:cav2} with the proximal
mixture $\pcmx{\gamma}(g_k,L_k)_{1\leq k\leq p}$ of \eqref{e:pcm2}.
Let us start with some basic inequalities.

\begin{proposition}
\label{p:2}
Let $\gamma\in\RPP$ and
\begin{equation}
\varphi\colon\HH\to\RX\colon x\mapsto
\inf_{\substack{
1\leq k\leq p\\ y_k^*\in\partial g_k(L_kx)}}
\brk4{\sum_{k=1}^p\alpha_k\qq_{\GG_k}(y_k^*)
-\qq_{\HH}\brk3{\sum_{k=1}^p\alpha_kL_k^*y_k^*}}.
\end{equation}
Then 
$\sum_{k=1}^p\alpha_k\brk{\lenv{g_k}{\gamma}}\circ L_k\leq
\pcmx{\gamma}(g_k,L_k)_{1\leq k\leq p}\leq
\cav(g_k,L_k)_{1\leq k\leq p}
\leq\pcmx{\gamma}(g_k,L_k)_{1\leq k\leq p}+\gamma\varphi$.
\end{proposition}
\begin{proof}
Define $\GG$, $L$, and $g$ as in \eqref{e:p0iii}, and recall from
Proposition~\ref{p:0}\ref{p:0iii} that
$\pcmx{\gamma}(g_k,L_k)_{1\leq k\leq p}=L\proxcc{\gamma}g$.
Further, $g\circ L=\cav(g_k,L_k)_{1\leq k\leq p}$, and
as seen after \eqref{e:prop1}, $\|L\|\leq 1$. It follows from
\cite[Proposition~3.20(ii)]{Eect26} that
\begin{equation}
\sum_{k=1}^p\alpha_k\brk{\lenv{g_k}{\gamma}}\circ L_k
=\brk{\lenv{g}{\gamma}}\circ L
\leq L\proxcc{\gamma}g\leq g\circ L,
\end{equation}
where $L\proxcc{\gamma}g=\pcmx{\gamma}(g_k,L_k)_{1\leq k\leq p}$
and $g\circ L=\cav(g_k,L_k)_{1\leq k\leq p}$. It remains to show
the rightmost inequality. Let $x\in\HH$. The result is clear if 
$x\notin\dom\varphi$. Now, suppose that $x\in\dom\varphi$ and let
$\boldsymbol{y}^*=(y_k^*)_{1\leq k\leq p}\in\bigtimes_{k=1}^p
\partial g_k(L_kx)$. Then
\eqref{e:sp}--\eqref{e:p0iii} yield 
$L^*\boldsymbol{y}^*=\sum_{k=1}^p\alpha_kL_k^*y_k^*$ and
$\qq_{\GG}(\boldsymbol{y}^*)=
\sum_{k=1}^p\alpha_k\qq_{\GG_k}(y_k^*)$.
On the other hand, by \cite[Proposition~16.9]{Livre1}, 
$\boldsymbol{y}^*\in
\bigtimes_{k=1}^p\partial g_k(L_kx)
=\partial g(Lx)$ and we therefore deduce from
\cite[Proposition~3.23(i)]{Eect26} that
\begin{equation}
0\leq\brk1{\cav(g_k,L_k)_{1\leq k\leq p}}(x)
-\brk1{\pcmx{\gamma}(g_k,L_k)_{1\leq k\leq p}}(x)\leq\gamma\brk1{
\qq_{\GG}(\boldsymbol{y}^*)-\qq_{\HH}(L^*\boldsymbol{y}^*)}.
\end{equation}
Hence, taking the infimum over $\boldsymbol{y}^*\in
\bigtimes_{k=1}^p\partial g_k(L_kx)$ yields the claim.
\end{proof}

The following result establishes asymptotic relations between 
standard composite averages and proximal comixtures.

\begin{theorem}
\label{t:4}
Let $(\gamma_n)_{n\in\NN}$ be a
sequence in $\RPP$ such that $\gamma_n\downarrow 0$
and let $x\in\HH$. For brevity, write
$\scav=\cav(g_k,L_k)_{1\leq k\leq p}$ and
$(\forall n\in\NN)$
$\spcm{\gamma_n}=\pcmx{\gamma_n}(g_k,L_k)_{1\leq k\leq p}$. 
Then the following hold:
\begin{enumerate}
\item
\label{t:4i}
$\spcm{\gamma_n}(x)\uparrow\scav(x)$.
\item
\label{t:4ii}
Let $(x_n)_{n\in\NN}$ be a sequence in $\HH$ such that
$x_n\weakly x$. Then
$\scav(x)\leq\varliminf\spcm{\gamma_n}(x_n)$.
\item
\label{t:4iii}
Suppose that $\scav$ is proper and let $\gamma\in\RPP$. Then
$\prox_{\gamma\spcm{\gamma_n}}x\to\prox_{\gamma\scav}x$.
\item
\label{t:4iv}
Let $(\varepsilon_n)_{n\in\NN}$ be a sequence in $\RPP$ such that
$\varepsilon_n\downarrow 0$, and suppose that there exists a
bounded sequence $(z_n)_{n\in\NN}$ such that
$(\forall n\in\NN)$ 
$z_n\in\varepsilon_n$--$\Argmin\brk{f+\spcm{\gamma_n}+h}$.
Then the following hold:
\begin{enumerate}
\item
\label{t:4iva}
$\inf_{x\in\HH}\brk{f(x)+\spcm{\gamma_n}(x)+h(x)}\uparrow
\min_{x\in\HH}\brk{f(x)+\scav(x)+h(x)}$.
\item
\label{t:4ivb}
Every weak sequential cluster point of $(z_n)_{n\in\NN}$ is in
$\Argmin\brk{f+\scav+\,h}$.
\end{enumerate}
\end{enumerate}
\end{theorem}
\begin{proof}
\ref{t:4i}: 
Define $L$ and $g$ as in \eqref{e:p0iii}, and recall from
Proposition~\ref{p:0}\ref{p:0iii} that $(\forall n\in\NN)$
$\spcm{\gamma_n}=L\proxcc{\gamma_n}g$. By items (ii) and (iv) in 
\cite[Theorem~3.29]{Eect26}, the function
$\NN\to\RX\colon n\mapsto\spcm{\gamma_n}(x)$ is
increasing and $\lim_{n\to\pinf}\spcm{\gamma_n}(x)=
g(Lx)=\scav(x)$.

\ref{t:4ii}: 
The weak continuity of bounded linear operators
\cite[Lemma~2.41]{Livre1} yields $(\forall k\in\{1,\ldots,p\})$
$L_kx_n\weakly L_kx$. Thus, \cite[Proposition~2.2(d)]{Vilches21}
implies that
\begin{equation}
\label{e:p4ii}
(\forall k\in\{1,\ldots,p\})\quad
g_k(L_kx)\leq\varliminf(\lenv{g_k}{\gamma_n})(L_kx_n).
\end{equation}
On the other hand, recall that 
$\scav=\cav(g_k,L_k)_{1\leq k\leq p}$. Therefore, it follows
from \eqref{e:p4ii}, \cite[Lemma~1.16]{Livre1}, and
Proposition~\ref{p:2} that
\begin{align}
\scav(x)\leq
\sum_{k=1}^p\alpha_k\varliminf(\lenv{g_k}{\gamma_n})(L_kx_n)
\leq\varliminf\sum_{k=1}^p\alpha_k(\lenv{g_k}{\gamma_n})(L_kx_n)
\leq\varliminf\spcm{\gamma_n}(x_n).
\end{align}

\ref{t:4iii}: 
Recall from Proposition~\ref{p:1}\ref{p:1i} that 
$(\forall n\in\NN)$ $\spcm{\gamma_n}\in\Gamma_0(\HH)$. Therefore,
the result follows from \ref{t:4i}, \ref{t:4ii},
the equivalence (i)$\Leftrightarrow$(iii) in
\cite[Proposition~3.19]{Atto84}, and
the equivalence (i)$\Leftrightarrow$(ii) in
\cite[Theorem~3.26]{Atto84}.

\ref{t:4iv}:
By Proposition~\ref{p:1}\ref{p:1i} and \ref{t:4i},
$(\spcm{\gamma_n})_{n\in\NN}$ is an increasing sequence of
functions in $\Gamma_0(\HH)$ with
$\sup_{n\in\NN}\spcm{\gamma_n}=\scav$. Therefore, the equivalence
(iii)$\Leftrightarrow$(iv) in \cite[Theorem~9.1]{Livre1} implies
that $(f+\spcm{\gamma_n}\!+h)_{n\in\NN}$ is an increasing
sequence of weakly lower semicontinuous convex functions with
supremum $f+\scav+h$. Since $(z_n)_{n\in\NN}$ is bounded,
\cite[Lemma~2.45]{Livre1} asserts that $(z_n)_{n\in\NN}$ 
admits a weakly convergent
subsequence. Hence, the results follow from
\cite[Proposition~2.42]{Atto84} applied to the weak topology.
\end{proof}

The following result focuses on the case in which the functions
$(g_k)_{1\leq k\leq p}$ are Lipschitzian.

\begin{proposition}
\label{p:3}
Assume that, for every $k\in\{1,\ldots,p\}$, $g_k$ is real-valued
and $\mu_k$-Lipschitzian for some $\mu_k\in\RPP$. Set
$\mu=\sum_{k=1}^p\alpha_k\mu_k\|L_k\|$, set
$\vartheta=(1/2)\sum_{k=1}^p\alpha_k\mu_k^2$, let $x\in\HH$, 
and let $\gamma\in\RPP$. Then the following hold:
\begin{enumerate}
\item
\label{p:3i}
$\pcmx{\gamma}(g_k,L_k)_{1\leq k\leq p}$ is $\mu$-Lipschitzian.
\item
\label{p:3ii}
$0\leq\cav(g_k,L_k)_{1\leq k\leq p}
-\pcmx{\gamma}(g_k,L_k)_{1\leq k\leq p}
\leq\gamma\vartheta$.
\item
\label{p:3iii}
There exists $r\in B(0;2\gamma\mu)$ such that
$\prox_{\gamma\pcmx{\gamma}(g_k,L_k)_{1\leq k\leq p}}x
=\prox_{\gamma\!\cav(g_k,L_k)_{1\leq k\leq p}}(x+r)$.
\item
\label{p:3v}
$\|\prox_{\gamma\!\cav(g_k,L_k)_{1\leq k\leq p}}x
-\prox_{\gamma\pcmx{\gamma}(g_k,L_k)_{1\leq k\leq p}}x\|
\leq\gamma\min\{2\mu,\sqrt{2\vartheta}\}$.
\item
\label{p:3vi}
Let $\sigma\in\RPP$ and assume that $h\in\Gamma_{-\sigma}(\HH)$.
Then
\begin{equation}
\label{e:45}
\left\|\argmin\brk1{\cav(g_k,L_k)_{1\leq k\leq p}+\,h}-\argmin
\brk1{\pcmx{\gamma}(g_k,L_k)_{1\leq k\leq p}+h}\right\|_{\HH}
\leq\dfrac{2\mu}{\sigma}.
\end{equation}
\end{enumerate}
\end{proposition}
\begin{proof}
According to \cite[Corollary~17.19]{Livre1}, a lower semicontinuous
convex function $\varphi\colon\HH\to\RR$ is $\mu$-Lipschitzian if
and only if $\ran\partial\varphi\subset B(0;\mu)$. Thus, the
Lipschitz continuity of the functions $(g_k)_{1\leq k\leq p}$
implies that
\begin{equation}
\label{e:15}
(\forall k\in\{1,\ldots,p\})\quad
\ran\partial g_k\subset B(0;\mu_k).
\end{equation}
Further, by \cite[Proposition~16.27]{Livre1}, 
$(\forall k\in\{1,\ldots,p\})$ $\dom\partial g_k=\GG_k$.

\ref{p:3i}:
By virtue of Proposition~\ref{p:1}\ref{p:1iiia+},
$\dom\pcmx{\gamma}(g_k,L_k)_{1\leq k\leq p}=\HH$ and
therefore \cite[Proposition~16.27]{Livre1} yields
$\dom\partial\,\pcmx{\gamma}(g_k,L_k)_{1\leq k\leq p}=\HH$. Now,
let $u\in\ran\partial\,\pcmx{\gamma}(g_k,L_k)_{1\leq k\leq p}$.
Then there exists $x\in\HH$ such that 
$u\in(\partial\,\pcmx{\gamma}(g_k,L_k)_{1\leq k\leq p})(x)$
which, by \eqref{e:jjm}, is equivalent to
$\prox_{\gamma\pcmx{\gamma}(g_k,L_k)_{1\leq k\leq p}}
(x+\gamma u)=x$.
Now, set $(\forall k\in\{1,\ldots,p\})$ 
$y_k^*=L_k(x+\gamma u)-\prox_{\gamma g_k}(L_k(x+\gamma u))$.
By Proposition~\ref{p:1}\ref{p:1ii},
\begin{align}
\label{e:15b}
x+\gamma u-\sum_{k=1}^p\alpha_kL_k^*y_k^*=x.
\end{align}
Further, it follows from \eqref{e:jjm} and \eqref{e:15} that
\begin{align}
\label{e:16}
(\forall k\in\{1,\ldots,p\})\quad\dfrac{1}{\gamma}y_k^*\in
\partial g_k\brk1{L_k(x+\gamma u)-y_k^*}\subset
\ran\partial g_k\subset B(0;\mu_k).
\end{align}
Therefore, it follows from \eqref{e:15b} and \eqref{e:16} that
\begin{align}
\|u\|_{\HH}
=\bigg\|\dfrac{1}{\gamma}
\sum_{k=1}^p\alpha_kL_k^*y_k^*\bigg\|_{\HH}
\leq\sum_{k=1}^p\alpha_k\|L_k\|
\left\|\dfrac{1}{\gamma}y_k^*\right\|_{\GG_k}
\leq\sum_{k=1}^p\alpha_k\|L_k\|\mu_k=\mu.
\end{align}
Altogether, 
$\ran\partial\,\pcmx{\gamma}(g_k,L_k)_{1\leq k\leq p}\subset
B(0;\mu)$, which yields the conclusion via
\cite[Corollary~17.19]{Livre1}. 

\ref{p:3ii}:
Recall that $(\forall k\in\{1,\ldots,p\}$) 
$\dom\partial g_k=\GG_k$. For every $k\in\{1,\ldots,p\}$, let 
$y_k^*\in\partial g_k(L_kx)$. Thus, it follows from
Proposition~\ref{p:2} and \eqref{e:15} that
\begin{equation}
0\leq\brk1{\cav(g_k,L_k)_{1\leq k\leq p}}(x)
-\brk1{\pcmx{\gamma}(g_k,L_k)_{1\leq k\leq p}}(x)
\leq\gamma\sum_{k=1}^p\alpha_k\qq_{\GG_k}(y_k^*) 
\leq\dfrac{\gamma}{2}\sum_{k=1}^p\alpha_k\mu_k^2
=\gamma\vartheta.
\end{equation}

\ref{p:3iii}:
Set $u_{\gamma}=
\prox_{\gamma\pcmx{\gamma}(g_k,L_k)_{1\leq k\leq p}}x$. 
Then, by \eqref{e:jjm}, 
$x-u_{\gamma}\in\gamma(\partial\,\pcmx{\gamma}
(g_k,L_k)_{1\leq k\leq p})(u_{\gamma})$. On the other hand, since
$\dom\partial\cav(g_k,L_k)_{1\leq k\leq p}=\HH$, there exists
$r\in\HH$ such that $x+r-u_{\gamma}\in
\gamma(\partial\cav(g_k,L_k)_{1\leq k\leq p})(u_{\gamma})$. Thus,
by \eqref{e:jjm}, $\prox_{\gamma\!\cav(g_k,L_k)_{1\leq k\leq
p}}(x+r)=u_{\gamma}$. Further, by \ref{p:3i} and
\cite[Corollary~17.19]{Livre1}, \begin{equation}
\ran\partial\cav(g_k,L_k)_{1\leq k\leq p}\subset B(0;\mu)
\quad\text{and}\quad\ran\partial\,
\pcmx{\gamma}(g_k,L_k)_{1\leq k\leq p}\subset B(0;\mu).
\end{equation}
Therefore, $r=(x+r-u_{\gamma})-(x-u_{\gamma})\in\gamma B(0;\mu)
+\gamma B(0;\mu)=B(0;2\gamma\mu)$.

\ref{p:3v}:
Set $u=\prox_{\gamma\cav(g_k,L_k)_{1\leq k\leq p}}x$, 
$u_{\gamma}=\prox_{\gamma\pcmx{\gamma}(g_k,L_k)_{1\leq k\leq p}}x$,
and
\begin{equation}
\label{e:p3v}
(\forall k\in\{1,\ldots,p\})\quad 
y_k^*=L_kx-\prox_{\gamma g_k}(L_kx).
\end{equation}
Further, set $z=\sum_{k=1}^p\alpha_kL_k^*y_k^*$.
Then, Proposition~\ref{p:1}\ref{p:1ii} yields $x-z=u_{\gamma}$,
whereas \eqref{e:jjm}, \eqref{e:p3v}, and \eqref{e:15} imply that
\begin{equation}
\label{e:40}
(\forall k\in\{1,\ldots,p\})\quad
y_k^*\in\gamma\partial g_k(L_kx-y_k^*)\subset B(0;\gamma\mu_k).
\end{equation}
By \cite[Corollaries~16.48(iii) and 16.53(i)]{Livre1},
$x-u\in\gamma(\partial\cav(g_k,L_k)_{1\leq k\leq p})(u)
=\gamma\sum_{k=1}^p\alpha_kL_k^*\brk{\partial g_k(L_ku)}$.
Thus, for every $k\in\{1,\ldots,p\}$, there exists 
$w_k^*\in\gamma\partial g_k(L_ku)$ such that 
$x-u=\sum_{k=1}^p\alpha_kL_k^*w_k^*$. It follows from \eqref{e:40}
and the monotonicity of the subdifferential operators 
$(\gamma\partial g_k)_{1\leq k\leq p}$ \cite[Theorem~20.25]{Livre1}
that
\begin{equation}
\label{e:41}
(\forall k\in\{1,\ldots,p\})\quad
\scal{L_kx-y_k^*-L_ku}{y_k^*-w_k^*}_{\GG_k}\geq 0.
\end{equation}
Since $\sum_{k=1}^p\alpha_k\|L_k\|^2\leq 1$ by
Assumption~\ref{a:1}, \eqref{e:40} implies that
\begin{align}
\label{e:42}
\sum_{k=1}^p\alpha_k\|y_k^*-L_kz\|_{\GG_k}^2
=\sum_{k=1}^p\alpha_k\|y_k^*\|_{\GG_k}^2
-2\|z\|_{\HH}^2+\sum_{k=1}^p\alpha_k\|L_kz\|_{\GG_k}^2
\leq\sum_{k=1}^p\alpha_k\|y_k^*\|_{\GG_k}^2
\leq2\gamma^2\vartheta,
\end{align}
whereas the Cauchy--Schwarz inequality yields
\begin{equation}
\label{e:43}
\|z\|_{\HH}^2\leq 
\brk3{\sum_{k=1}^p\alpha_k\|L_k\|\|y_k^*\|_{\GG_k}}^2
\leq\brk3{\sum_{k=1}^p\alpha_k\|L_k\|^2}
\brk3{\sum_{k=1}^p\alpha_k\|y_k^*\|_{\GG_k}^2}
\leq\sum_{k=1}^p\alpha_k\|y_k^*\|_{\GG_k}^2.
\end{equation}
Altogether, we deduce from $x=z+u_{\gamma}$, 
$\sum_{k=1}^p\alpha_kL_k^*y_k^*=z$, 
$\sum_{k=1}^p\alpha_kL_k^*w_k^*=x-u$, \eqref{e:41}, the
Cauchy--Schwarz inequality, \eqref{e:42}, and \eqref{e:43} that
\begin{align}
\label{e:50}
0&\leq\sum_{k=1}^p\alpha_k\scal{L_ku_{\gamma}-L_ku+L_kz-y_k^*}{
y_k^*-w_k^*}_{\GG_k}\nonumber\\
&=\scal{u_{\gamma}-u}{(x-u_{\gamma})-(x-u)}_{\HH}
+\sum_{k=1}^p\alpha_k\scal{L_kz-y_k^*}{y_k^*-w_k^*}_{\GG_k}
\nonumber\\
&=-\|u-u_{\gamma}\|_{\HH}^2
+\sum_{k=1}^p\alpha_k\scal{y_k^*-L_kz}{w_k^*}_{\GG_k}
-\sum_{k=1}^p\alpha_k\scal{y_k^*-L_kz}{y_k^*}_{\GG_k}\nonumber\\
&\leq-\|u-u_{\gamma}\|_{\HH}^2
+\sqrt{\sum_{k=1}^p\alpha_k\|y_k^*-L_kz\|_{\GG_k}^2}
\,\sqrt{\sum_{k=1}^p\alpha_k\|w_k^*\|_{\GG_k}^2}
+\brk3{\|z\|_{\HH}^2-\sum_{k=1}^p\alpha_k\|y_k^*\|_{\GG_k}^2}
\nonumber\\
&\leq-\|u-u_{\gamma}\|_{\HH}^2+\sqrt{2\gamma^2\vartheta}
\sqrt{\gamma^2\sum_{k=1}^p\alpha_k\mu_k^2}\nonumber\\
&=-\|u-u_{\gamma}\|_{\HH}^2+2\gamma^2\vartheta.
\end{align}
However, by \ref{p:3iii}, there exists 
$r\in B(0;2\gamma\mu)$ such that
$\prox_{\gamma\pcmx{\gamma}(g_k,L_k)_{1\leq k\leq p}}x
=\prox_{\gamma\!\cav(g_k,L_k)_{1\leq k\leq p}}(x+r)$.
Hence, the nonexpansiveness of the proximity operator
\cite[Proposition~12.28]{Livre1} implies that
\begin{equation}
\label{e:51}
\|u-u_{\gamma}\|_{\HH}=\bigl\|\,\prox_{\gamma\!
\cav(g_k,L_k)_{1\leq k\leq p}}x
-\prox_{\gamma\!\cav(g_k,L_k)_{1\leq k\leq p}}(x+r)\bigr\|_{\HH}
\leq\|r\|_{\HH}\leq 2\gamma\mu.
\end{equation}
Finally, the result follows from \eqref{e:50} and
\eqref{e:51}.

\ref{p:3vi}:
Note that $\cav(g_k,L_k)_{1\leq k\leq p}+h\in\Gamma_{-\sigma}(\HH)$ 
and, by Proposition~\ref{p:1}\ref{p:1i}, also
$\pcmx{\gamma}(g_k,L_k)_{1\leq k\leq p}+h\in\Gamma_{-\sigma}(\HH)$.
Thus, \cite[Corollary~11.17]{Livre1} implies that
$\cav(g_k,L_k)_{1\leq k\leq p}+h$ and
$\pcmx{\gamma}(g_k,L_k)_{1\leq k\leq p}+h$ admit unique
minimizers.
Set $x=\argmin(\cav(g_k,L_k)_{1\leq k\leq p}+h)$ and
$x_\gamma=\argmin(\pcmx{\gamma}(g_k,L_k)_{1\leq k\leq p}+h)$.
Then the equivalence (i)$\Leftrightarrow$(viii) in 
\cite[Corollary~27.3]{Livre1} yields
\begin{equation}
\label{e:c45i}
x=\prox_{\gamma\!\cav(g_k,L_k)_{1\leq k\leq p}}
\brk1{x-\gamma\nabla h(x)}\quad\text{and}\quad x_{\gamma}
=\prox_{\gamma\pcmx{\gamma}(g_k,L_k)_{1\leq k\leq p}}\brk1{
x_{\gamma}-\gamma\nabla h(x_{\gamma})}.
\end{equation}
On the other hand, \ref{p:3iii} asserts that there exists $r\in
B(0;2\gamma\mu)$ such that 
\begin{equation}
\label{e:c45ii}
\prox_{\gamma\pcmx{\gamma}(g_k,L_k)_{1\leq k\leq p}}\brk1{
x-\gamma\nabla h(x)}
=\prox_{\gamma\!\cav(g_k,L_k)_{1\leq k\leq p}}\brk1{
x-\gamma\nabla h(x)+r}.
\end{equation}
Altogether, \eqref{e:c45i}, \eqref{e:c45ii}, and
the firm nonexpansiveness of 
$\prox_{\gamma\!\cav(g_k,L_k)_{1\leq k\leq p}}$
\cite[Proposition~12.28]{Livre1} yield
\begin{align}
\|x-x_{\gamma}\|_{\HH}^2
&=\left\|\prox_{\gamma\!\cav(g_k,L_k)_{1\leq k\leq p}}\brk1{
x-\gamma\nabla h(x)}
-\prox_{\gamma\pcmx{\gamma}(g_k,L_k)_{1\leq k\leq p}}\brk1{
x_{\gamma}-\gamma\nabla h(x_{\gamma}) }\right\|_{\HH}^2\nonumber\\
&=\left\|\prox_{\gamma\!\cav(g_k,L_k)_{1\leq k\leq p}}\brk1{
x-\gamma\nabla h(x)}
-\prox_{\gamma\!\cav(g_k,L_k)_{1\leq k\leq p}}\brk1{
x_{\gamma}-\gamma\nabla h(x_{\gamma})+r}\right\|_{\HH}^2\nonumber\\
&\leq\scal{x-\gamma\nabla h(x)-x_{\gamma}+\gamma\nabla
h(x_{\gamma})-r}{x-x_{\gamma}}_{\HH}\nonumber\\
&=\|x-x_{\gamma}\|_{\HH}^2
-\gamma\scal{\nabla h(x)-\nabla h(x_{\gamma})}{x-x_{\gamma}}_{\HH}
-\scal{r}{x-x_{\gamma}}_{\HH}.
\end{align}
Therefore, it follows from \cite[Example~22.4(iv)]{Livre1}, the
Cauchy--Schwarz inequality, and the inequality
$\|r\|_{\HH}\leq 2\gamma\mu$ that
\begin{align}
\gamma\sigma\|x-x_{\gamma}\|_{\HH}^2\leq
\gamma\scal{\nabla h(x)-\nabla h(x_{\gamma})}{x-x_{\gamma}}_{\HH}
\leq 2\gamma\mu\|x-x_{\gamma}\|_{\HH},
\end{align}
which implies \eqref{e:45}.
\end{proof}

As a special case of the above results, we recover some properties
of the proximal average \cite{Baus08,Kami21,Yuyl13}.

\begin{remark}[proximal average]
\label{r:17}
Suppose that $\sum_{k=1}^p\alpha_k=1$ and that, for every 
$k\in\{1,\ldots,p\}$, $\GG_k=\HH$ and $L_k=\Id_{\HH}$.
Then $\pcmx{\gamma}(g_k,L_k)_{1\leq k\leq p}
=\pave{\gamma}(g_k)_{1\leq k\leq p}$ is the proximal average of
\eqref{e:pav}, while 
$\cav(g_k,L_k)_{1\leq k\leq p}=\ave(g_k)_{1\leq k\leq p}$ is the
standard average of \eqref{e:avg}.
In this context, we recover the following results:
\begin{enumerate}
\item
Proposition~\ref{p:1}\ref{p:1i} yields
$\pave{\gamma}(g_k)_{1\leq k\leq p}\in\Gamma_0(\HH)$
(see \cite[Corollary~5.2]{Baus08}).
\item
Proposition~\ref{p:1}\ref{p:1ii} yields
$\prox_{\gamma\pave{\gamma}(g_k)_{1\leq k\leq p}}
=\sum_{k=1}^p\alpha_k\prox_{\gamma g_k}$
(see \cite[Theorem~6.7]{Baus08}).
\item
Proposition~\ref{p:1}\ref{p:1iv} yields
$\lenv{\pave{\gamma}(g_k)_{1\leq k\leq p}}{\gamma}
=\sum_{k=1}^p\alpha_k\lenv{g_k}{\gamma}$ 
(see \cite[Theorem~6.2(i)]{Baus08}).
\item
Proposition~\ref{p:2} yields
$\pave{\gamma}(g_k)_{1\leq k\leq p}\leq\ave$
(see \cite[Theorem~5.4]{Baus08}).
\item
Let $(\gamma_n)_{n\in\NN}$ be a sequence in $\RPP$ such that
$\gamma_n\downarrow 0$. 
\begin{enumerate}
\item
Theorem~\ref{t:4}\ref{t:4i} yields
$\pave{\gamma_n}(g_k)_{1\leq k\leq p}\uparrow
\ave(g_k)_{1\leq k\leq p}$ (see \cite[Theorem~8.5]{Baus08}).
\item
Theorem~\ref{t:4}\ref{t:4i}--\ref{t:4ii} imply that
$(\pave{\gamma_n}(g_k)_{1\leq k\leq p})_{n\in\NN}$
epi-converges to $\ave(g_k)_{1\leq k\leq p}$
(see \cite[Corollary~9.6]{Baus08}).
\end{enumerate}
\item
Suppose that, for every $k\in\{1,\ldots,p\}$, $g_k$ is real-valued
and $\mu_k$-Lipschitzian for some $\mu_k\in\RPP$. Set
$\mu=\sum_{k=1}^p\alpha_k\mu_k$ and 
$\vartheta=(1/2)\sum_{k=1}^p\alpha_k\mu_k^2$.
\begin{enumerate}
\item
Proposition~\ref{p:3}\ref{p:3ii} yields
$0\leq\ave(g_k)_{1\leq k\leq p}-\pave{\gamma}(g_k)_{1\leq k\leq p}
\leq\gamma\vartheta$
(see \cite[Proposition~4]{Yuyl13}).
\item
Proposition~\ref{p:3}\ref{p:3v} yields
$\|\prox_{\gamma\pave{\gamma}(g_k)_{1\leq k\leq p}}
-\prox_{\gamma\!\ave(g_k)_{1\leq k\leq p}}\|_{\HH}\leq 2\gamma\mu$
(see \cite[p.~851]{Kami21}).
\end{enumerate}
\end{enumerate}
In addition, we obtain the following new properties of the proximal
average.
\begin{enumerate}
\setcounter{enumi}{6}
\item
Let $(\gamma_n)_{n\in\NN}$ be a sequence in $\RPP$ such that
$\gamma_n\downarrow 0$. It follows from
Theorem~\ref{t:4}\ref{t:4iii} that
$\prox_{\gamma\pave{\gamma_n}(g_k)_{1\leq k\leq p}}x
\to\prox_{\gamma\!\ave(g_k)_{1\leq k\leq p}}x$.
\item
Suppose that, for every $k\in\{1,\ldots,p\}$, $g_k$ is real-valued
and $\mu_k$-Lipschitzian for some $\mu_k\in\RPP$. Set
$\mu=\sum_{k=1}^p\alpha_k\mu_k$. Then, by
Proposition~\ref{p:3}\ref{p:3iii}, there exists
$r\in B(0;2\gamma\mu)$ such that
$\prox_{\gamma\pave{\gamma}(g_k)_{1\leq k\leq p}}x
=\prox_{\gamma\!\ave(g_k)_{1\leq k\leq p}}(x+r)$.
\end{enumerate}
\end{remark}

Generalization of Problems~\ref{prob:1} and \ref{prob:3} to
arbitrary families of functions and linear operators can be
formulated as follows.

\begin{remark}[integral proximal comixture]
Let $(\Omega,\FF,\mu)$ be a complete $\sigma$-finite measure space
and let $\HS$ be a separable real Hilbert space. For every 
$\omega\in\Omega$, let $\GW$ be a real Hilbert space, let
$\gw\in\Gamma_0(\GW)$, and let $\LW\colon\HS\to\GW$ be a bounded
linear operator. Assume that
$0<\int_{\Omega}\norm{\LW}^2\mu(d\omega)\leq 1$.
Under mild assumptions, the \emph{integral proximal comixture}
of $(\gw)_{\omega\in\Omega}$ and
$(\LW)_{\omega\in\Omega}$ with parameter
$\gamma\in\RPP$ is \cite[Definition~4.2]{Jota24} 
\begin{equation}
\label{e:ipm}
\Rcm{\gamma}(\gw,\LW)_{\omega\in\Omega}
=\brk2{\mathsf{\varphi}^*-\dfrac{\gamma}{2}\|\cdot\|_{\HS}^2}^*,
\quad\text{where}\quad(\forall\mathsf{x}\in\HS)\quad
\mathsf{\varphi}(\mathsf{x})=\int_{\Omega}
(\lenv{\gw}{\gamma})(\LW\mathsf{x})\mu(d\omega).
\end{equation}
Now let $\mathsf{f}\in\Gamma_0(\HS)$ and let
$\mathsf{h}\colon\HS\to\RR$ be a convex and differentiable 
function with Lipschitzian gradient. Then a generalization of
Problem~\ref{prob:1} is
\begin{equation}
\label{e:p6}
\minimize{\mathsf{x}\in\HS}{\mathsf{f}(\mathsf{x})
+\int_{\Omega}\gw(\LW\mathsf{x})\,\mu(d\omega)
+\mathsf{h}(\mathsf{x})}
\end{equation}
and a generalization of Problem~\ref{prob:3} is
\begin{equation}
\label{e:p7}
\minimize{\mathsf{x}\in\HS}{\mathsf{f}(\mathsf{x})
+\brk2{\Rcm{\gamma}(\gw,\LW)_{\omega\in\Omega}}(\mathsf{x})
+\mathsf{h}(\mathsf{x})}.
\end{equation}
Indeed, we recover \eqref{e:p1} from \eqref{e:p6} 
and \eqref{e:p3} from \eqref{e:p7} by taking
$\Omega=\{1,\ldots,p\}$ and $\FF=2^{\Omega}$, and setting 
$(\forall k\in\{1,\ldots,p\})$ $\mu(\{k\})=\alpha_k$.
Most of the results of this section extend to this abstract
framework using the tools of \cite{Jota24,Eect26}. If $\mu$ is a
probability measure, we can regard \eqref{e:p1} and \eqref{e:p3} as
empirical versions of \eqref{e:p6} and \eqref{e:p7}, respectively.
\end{remark}

\section{Algorithms}
\label{sec:4}

Several splitting methods are available to solve
Problems~\ref{prob:1} and \ref{prob:3} \cite{Acnu24}. For solving
Problem~\ref{prob:1}, methods that require the proximity operator
of $\cav(g_k,L_k)_{1\leq k\leq p}$ are ruled out as this operator
is not tractable. We have therefore considered those of
\cite{Acnu24,Svva12,Cond13,Bang13,Yan18} and found that, overall,
the following method \cite{Yan18} (see also \cite{Cond22}), which
is based on \cite{Davi17}, performed the best in our experiments.

\begin{proposition}
\label{p:ave3}
Consider the setting of Problem~\ref{prob:1} and assume that
\begin{equation}
0\in\ran\brk2{\partial f+\sum_{k=1}^p\alpha_kL_k^*\circ\brk{
\partial g_k}\circ L_k+\nabla h}. 
\end{equation} 
Let $\eta\in\intv[o]{0}{2\beta}$ and
$\sigma\in\intv[o]{0}{1/(\eta\sum_{k=1}^p\alpha_k\|L_k\|^2)}$. 
Let $q_0\in\HH$ and, for every $k\in\{1,\ldots,p\}$, let
$u_{k,0}^*\in\GG_k$. Iterate
\begin{equation}
\label{e:12}
\hskip -2.7mm
\begin{array}{l}
\text{for}\;n=0,1,\ldots\\
\left\lfloor
\begin{array}{l}
x_n=\prox_{\eta f}\brk2{\eta\brk1{q_n-\sum_{k=1}^p\alpha_k
L_k^*u_{k,n}^*}}\\
q_{n+1}=\dfrac{1}{\eta}x_n-\nabla h(x_n)\\
\text{for}\;k=1,\ldots,p\\
\left\lfloor
\begin{array}{l}
r_{k,n}=u_{k,n}^*+\sigma L_k\brk1{x_n+\eta(q_{n+1}-q_n)}\\
u_{k,n+1}^*=r_{k,n}-\sigma\prox_{g_k/\sigma}(r_{k,n}/\sigma).
\end{array} 
\right.\\ 
\end{array} 
\right. 
\end{array}
\end{equation} 
Then $(x_n)_{n\in\NN}$ converges weakly to a solution to 
Problem~\ref{prob:1}. 
\end{proposition}
\begin{proof}
Apply \cite[Theorem~1]{Cond22} in the Hilbert space $\GG$ of
\eqref{e:sp}, with $g$ and $L$ defined as in \eqref{e:p0iii}.
\end{proof}

We now turn to Problem~\ref{prob:3}, for which the method proposed
in \cite[Section~3.1]{Davi17} can be implemented directly.
Given $y_0\in\HH$, this algorithm iterates
\begin{equation}
\label{e:11}
\hskip -2.7mm
\begin{array}{l}
\text{for}\;n=0,1,\ldots\\
\left\lfloor
\begin{array}{l}
x_n=\prox_{\gamma\pcmx{\gamma}(g_k,L_k)_{1\leq k\leq p}}y_n\\
z_n=\prox_{\gamma f}\brk1{2x_n-y_n-\gamma\nabla h(x_n)}\\
y_{n+1}=y_n+\lambda_n(z_n-x_n).
\end{array}
\right.
\end{array}
\end{equation}
Hence, we arrive at the following method.

\begin{proposition}
\label{p:pcm}
Consider the setting of Problem~\ref{prob:3}. Assume that
$\gamma\in\intv[o]{0}{2\beta}$ and that
\begin{equation}
\label{e:ppcm}
0\in\ran\brk1{
\partial f+\partial\,
\pcmx{\gamma}(g_k,L_k)_{1\leq k\leq p}+\nabla h}.
\end{equation}
Set $\delta=2-\gamma/(2\beta)$, let
$(\lambda_n)_{n\in\NN}$ be a sequence in $\intv[o]{0}{\delta}$ such
that $\sum_{n\in\NN}\lambda_n(\delta-\lambda_n)=\pinf$, and let
$y_0\in\HH$. Iterate
\begin{equation}
\label{e:13}
\hskip -2.7mm
\begin{array}{l}
\text{for}\;n=0,1,\ldots\\
\left\lfloor
\begin{array}{l}
x_n=y_n-\sum_{k=1}^p\alpha_kL_k^*
\brk1{L_ky_n-\prox_{\gamma g_k}(L_ky_n)}\\
z_n=\prox_{\gamma f}\brk1{2x_n-y_n-\gamma\nabla h(x_n)}\\
y_{n+1}=y_n+\lambda_n(z_n-x_n).
\end{array}
\right.
\end{array}
\end{equation}
Then $(x_n)_{n\in\NN}$ converges weakly to a solution to
Problem~\ref{prob:3}.
\end{proposition}
\begin{proof}
By Proposition~\ref{p:1}\ref{p:1ii}, the algorithms \eqref{e:11}
and \eqref{e:13} are equivalent. Therefore, the assertion follows
from \cite[Theorem~2.1.1(b)]{Davi17} applied to \eqref{e:11}.
\end{proof}

Heuristically, one can expect \eqref{e:13} to yield fast
convergence than \eqref{e:12}. Indeed, since the proximity operator
of the aggregation $\pcmx{\gamma}(g_k,L_k)_{1\leq k\leq p}$ is
computable in closed form via Proposition~\ref{p:1}\ref{p:1ii},
algorithm~\eqref{e:11} processes only the three functions 
$f$, $\pcmx{\gamma}(g_k,L_k)_{1\leq k\leq p}$, and $h$. It also
requires minimum variable storage. By contrast, the standard
composite average function $\cav(g_k,L_k)_{1\leq k\leq p}$ in
Problem~\ref{prob:1} has no explicit proximity operator and its 
processing in algorithm~\eqref{e:12} requires the splitting of the
$p$ functions $(g_k)_{1\leq k\leq p}$ and the $p$ linear operator
$(L_k)_{1\leq k\leq p}$ (the same holds true for all splitting
methods involving standard composite averages \cite{Acnu24}).
Additionally, algorithm~\eqref{e:12} requires the storage of
significantly more variables than \eqref{e:13}, in particular of
the dual variables $(u^*_{k,n})_{1\leq k\leq p}$.

\section{Numerical experiments}
\label{sec:5}

We illustrate the proposed proximal comixture model of
Problem~\ref{prob:3} through several experiments. The algorithms
are executed with all initial vectors set to $0$ and they use
parameters $\eta$, $\sigma$, and $(\lambda_n)_{n\in\NN}$ aiming at
optimizing their performance. The stopping criterion to approximate
the limit $x_{\infty}$ of the sequence $(x_n)_{n\in\NN}$ produced
by an algorithm is $\|x_{n+1}-x_n\|/\|x_n\|<10^{-6}$. The
simulations are run using GNU Octave $8.4.0$ on a laptop running
Ubuntu $24.04.4$. The code is available on the GitHub repository
\href{https://github.com/djcornejo/pcmm}
{https://github.com/djcornejo/pcmm}.
The results of
the these experiments support the fact that the proximal comixture
model of Problem~\ref{prob:3} leads to reliable solutions which can
be obtained by efficient algorithms that are much less demanding 
in terms of memory requirements than those of Problem~\ref{prob:1}.

\subsection{Experiment 1: Proximal total variation denoising}
\label{sec:52}

The purpose of this experiment is to illustrate the asymptotic
property of Theorem~\ref{t:4}\ref{t:4iii} on a simple denoising
problem. Let $\overline{x}\in\RR^N$ $(N=256)$ be the original
$1$-dimensional signal shown in Fig.~\ref{fig:ex4im1}(a) and let
\begin{equation}
\label{e:d8}
z=\overline{x}+\dfrac{1}{2}w
\end{equation}
be the noisy observation shown in Fig.~\ref{fig:ex4im1}(b), where
$w\in\RR^N$ is a realization of a normalized white Gaussian noise 
vector. We denote by 
$D\colon\RR^N\to\RR^N\colon(\xi_1,\ldots,\xi_N)\mapsto
(\xi_2-\xi_1,\ldots,\xi_N-\xi_{N-1},\xi_1-\xi_N)/2$ the normalized
discrete gradient operator.

\begin{figure}[h!]
\centering
\begin{tikzpicture}[scale=1]
\begin{axis}[height=0.23\textwidth,width=\textwidth,legend
style={at={(0.02,0.02)},legend cell align={left},
anchor=south west}, xmin =0, xmax=255, xtick distance = 20, 
ymin=-3, ymax=6, ytick={-3,-2,-1,0,1,2,3,4,5,6}, 
yticklabels={,-2,,0,,2,,4,,6}, grid=both]
\addplot [line width=0.5mm, mark=none, color=dblue] 
table[x={ind}, y={value}]
{figures/ptv/ex1/signal.txt};
\end{axis}
\end{tikzpicture}\\
(a) \\ 
\begin{tikzpicture}[scale=1]
\begin{axis}[height=0.23\textwidth,width=\textwidth,legend
style={at={(0.02,0.02)},legend cell align={left},
anchor=south west}, xmin =0, xmax=255, xtick distance = 20, 
ymin=-3, ymax=7, ytick={-4,-3,-2,-1,0,1,2,3,4,5,6},
yticklabels={-4,,-2,,0,,2,,4,,6}, grid=both]
\addplot [line width=0.5mm, mark=none, color=dblue] 
table[x={ind}, y={value}]
{figures/ptv/ex1/noisy.txt};
\end{axis}
\end{tikzpicture}\\
(b) \\
\begin{tikzpicture}[scale=1]
\begin{axis}[height=0.23\textwidth,width=\textwidth,legend
style={at={(0.02,0.02)},legend cell align={left},
anchor=south west}, xmin =0, xmax=255, xtick distance = 20, 
ymin=-3, ymax=6, ytick={-4,-3,-2,-1,0,1,2,3,4,5,6},
yticklabels={-4,,-2,,0,,2,,4,,6}, grid=both]
\addplot [line width=0.5mm, mark=none, color=dgreen] 
table[x={ind}, y={value}]
{figures/ptv/ex1/x_ptv.txt};
\end{axis}
\end{tikzpicture}\\
(c) \\ [-0.5em]
\caption{
(a) Original signal $\overline{x}$. 
(b) Noisy observation $z$. 
(c) Solution to Problem~\ref{p:ex4b} for $\gamma=10^{-3}$ 
(the solution to Problem~\ref{p:ex4a} is essentially identical).
}
\label{fig:ex4im1}
\end{figure}

\begin{problem}
\label{p:ex4a}
Let $\rho=3$ and let $\text{tv}=\|\cdot\|_1\circ D$ be the
standard total variation loss. The task is to
\begin{equation}
\label{e:ex4a}
\minimize{x\in\RR^N}{\text{tv}(x)+\dfrac{1}{2\rho}\|x-z\|^2}.
\end{equation}
\end{problem}

\begin{problem}
\label{p:ex4b}
Let $\rho=3$ and $\gamma\in\RPP$. Let us
introduce the \emph{proximal total variation} loss
$\text{ptv}_{\gamma}={D}\proxcc{\gamma}\|\cdot\|_1$. 
The task is to
\begin{equation}
\label{e:ex4b}
\minimize{x\in\RR^N}{\text{ptv}_{\gamma}(x)
+\dfrac{1}{2\rho}\|x-z\|^2}.
\end{equation}
\end{problem}

Problems~\ref{p:ex4a} and \ref{p:ex4b} are particular instances of 
Problems~\ref{prob:1} and \ref{prob:3}, respectively, where
$\HH=\RR^N$, $f=0$, $h=\|\cdot-z\|^2/(2\rho)$, $p=1$,
$\alpha_1=1$, $g_1=\|\cdot\|_1$, and $L_1=D$. In view of
Proposition~\ref{p:1}\ref{p:1i} and \eqref{e:prox}, the unique
solutions to Problems~\ref{p:ex4a} and \ref{p:ex4b} are,
respectively, $\prox_{\rho\,\text{tv}}z$ and
$\prox_{\rho\,\text{ptv}_{\gamma}}z$. Therefore,
Theorem~\ref{t:4}\ref{t:4iii} asserts that the solution curve
$(\prox_{\rho\,\text{ptv}_{\gamma}}z)_{\gamma\in\RPP}$ in
Problem~\ref{p:ex4b} converges to the solution to
Problem~\ref{p:ex4a}, to wit,
\begin{equation}
\label{e:ex4}
\prox_{\rho\,\text{ptv}_{\gamma}}z\to
\prox_{\rho\,\text{tv}}z\quad\text{as}\quad\gamma\downarrow 0.
\end{equation}
Since the proximity operator of \text{tv} is not known
explicitly, algorithm \eqref{e:12} can be applied to obtain the
solution to Problem~\ref{p:ex4a} On the other hand, algorithm
\eqref{e:13} can be used to obtain the solution to
Problem~\ref{p:ex4b}. Our numerical experiments confirmed that, for
$\gamma\leq 10^{-3}$, the solutions to Problems~\ref{p:ex4a} and
\ref{p:ex4b} were essentially identical, with a relative
difference
$\|\prox_{\rho\,\text{tv}}z-\prox_{\rho\,\text{ptv}_{10^{-3}}}z\|
/\|\prox_{\rho\,\text{tv}}z\|\approx 1.726\times 10^{-3}$, 
in conformity with \eqref{e:ex4}; see Fig.~\ref{fig:ex4im1}(c) for
the denoised signal.

\subsection{Experiment 2: Multiview image reconstruction}
\label{sec:ex3}

\begin{figure}[t]
\centering
\setlength{\tabcolsep}{0.2pt} 
\noindent
\begin{tabular}{ccc}
\includegraphics[width=4.6cm,height=4.5cm,clip]{
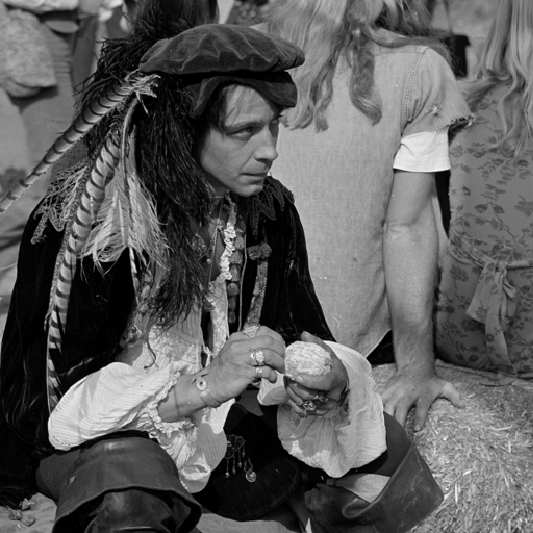} 
&\includegraphics[width=4.6cm,height=4.5cm,clip]{
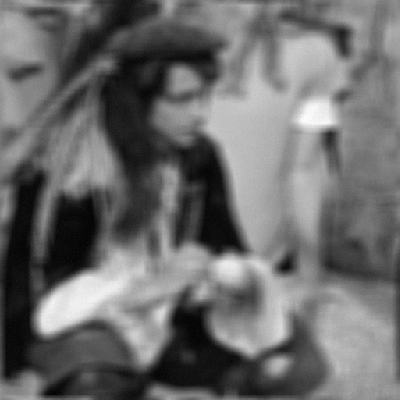} 
&\includegraphics[width=4.6cm,height=4.5cm,clip]{
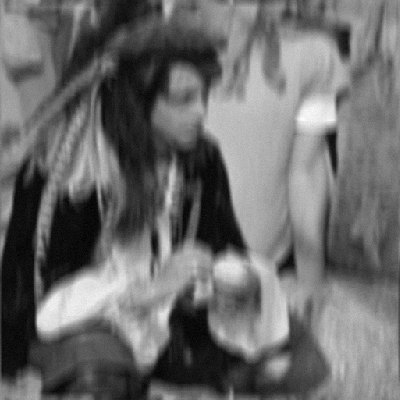} 
\\
\small{(a)} & \small{(b)} & \small{(c)}
\end{tabular} 
\vskip -0.2cm
\caption{
(a) Original image $\overline{x}$. 
(b) Degraded image $z_1$.
(c) Degraded image $z_2$.}
\label{fig:ex3im1}
\end{figure}

We address the problem of reconstructing the original image 
$\overline{x}\in C=[0,255]^N$ ($N=512^2$) of
Fig.~\ref{fig:ex3im1}(a) from a partial observation of its
possibly corrupted diffraction $r\approx\widehat{\overline{x}}$
over some frequency range $R$ \cite{Seza83}, where 
$\widehat{\overline{x}}$ denotes the two-dimensional discrete
Fourier transform of $\overline{x}$. To exploit this imprecise
information, we use the soft constraint distance penalty $d_E$,
where
\begin{equation}
\label{e:E}
E=\menge{x\in\RR^N}{(\forall\nu\in R)\;\;\widehat{x}(\nu)=r(\nu)}.
\end{equation}
The set $R$ contains the frequencies in 
$\{0,\ldots,15\}^2$ as well as those resulting from the 
symmetry properties of the discrete Fourier transform. Further,
two blurred noisy observations of $\overline{x}$ are
available, namely (see Fig.~\ref{fig:ex3im1}(b)--(c))
\begin{equation}
z_1=H_1\overline{x}+w_1
\quad\text{and}\quad 
z_2=H_2\overline{x}+w_2,
\end{equation}
where $H_1$ and $H_2$ model convolutional blurs with uniform
rectangular kernels of sizes $14\times 18$ and $20\times 5$,
respectively, while $w_1$ and $w_2$ represent zero-mean
white Gaussian noise
realizations of standard deviations of $2$ and $3$, respectively.
The blurred image-to-noise ratios are 34.51 dB and 31.06 dB,
respectively. Let
\begin{equation}
\label{e:dgrad}
D\colon\RR^N\to\RR^N\times\RR^N\colon x\mapsto(D_1x,D_2x),
\end{equation}
where $D_1$ and $D_2$ denote, respectively, the horizontal and
vertical first order discrete difference operators, i.e.,
\begin{equation}
\label{e:dgrad2}
(D_1x)_{i,j}=\begin{cases}
x_{i,j+1}-x_{i,j},&\text{if}\,\,j<512;\\
0,&\text{if}\,\,j=512,
\end{cases}
\quad\text{and}\quad
(D_2x)_{i,j}=\begin{cases}
x_{i+1,j}-x_{i,j},&\text{if}\,\,i<512;\\
0,&\text{if}\,\,i=512.
\end{cases}
\end{equation}
Let $\|\cdot\|_{1,2}\colon\RR^N\times\RR^N\to\RR\colon
(\xi_i,\eta_i)_{1\leq i\leq N}\mapsto
\sum_{i=1}^N\sqrt{\xi_i^2+\eta_i^2}$ and let
\begin{equation}
\mathfrak{h}_{\rho}\colon\RR\to\RR\colon\xi\mapsto
\begin{cases}
\rho\lvert\xi\rvert-\dfrac{\rho^2}{2},&\text{if}\quad
\lvert\xi\rvert>\rho;\\
\dfrac{\lvert\xi\rvert^2}{2},&\text{if}\quad
\lvert\xi\rvert\leq\rho
\end{cases}
\end{equation}
be the Huber function with parameter $\rho\in\RPP$.
Our first formulation involves a standard composite average.

\begin{problem}
\label{p:ex3a}
Let $\rho_1=3000$ and $\rho_2=4000$. The task is to 
\begin{equation}
\label{e:ex3a}
\minimize{x\in C}{
\brk2{\dfrac{1}{2}d_E(x)+\dfrac{1}{2}\|Dx\|_{1,2}}
+\brk2{\mathfrak{h}_{\rho_1}\brk1{\|H_1x-z_1\|}
+\mathfrak{h}_{\rho_2}\brk1{\|H_2x-z_2\|}}}.
\end{equation}
\end{problem}

The second formulation involves a proximal comixture.

\begin{problem}
\label{p:ex3c}
Let $\rho_1=3000$, $\rho_2=4000$, and $\gamma\in\intv[o]{0}{1}$.
The task is to 
\begin{equation}
\label{e:ex3c}
\minimize{x\in C}{\brk2{\pcmx{\gamma}\brk1{d_E,\Id;
\sqrt{8}\|\cdot\|_{1,2},D/\sqrt{8}}}(x)
+\brk2{\mathfrak{h}_{\rho_1}\brk1{\|H_1x-z_1\|}
+\mathfrak{h}_{\rho_2}\brk1{\|H_2x-z_2\|}}}.
\end{equation}
\end{problem}

Problems~\ref{p:ex3a} and \ref{p:ex3c} are particular instances of
Problems~\ref{prob:1} and \ref{prob:3}, respectively, where
$\HH=\RR^N$, $f=\iota_C$,
$h=\mathfrak{h}_{\rho_1}\circ\|H_1\cdot-z_1\|
+\mathfrak{h}_{\rho_2}\circ\|H_2\cdot-z_2\|$, $\beta=1/2$,
$p=2$, $\alpha_1=\alpha_2=1/2$, $\GG_1=\RR^N$, $g_1=d_E$,
$L_1=\Id$, $\GG_2=\RR^N\times\RR^N$,
$g_2=\sqrt{8}\|\cdot\|_{1,2}$, and $L_2=D/\sqrt{8}$.
In this case, $\|L_1\|^2=\|L_2\|^2=1$ and
\begin{equation}
\label{e:pbox}
\prox_{\gamma f}\colon\RR^N\to\RR^N\colon(\xi_i)_{1\leq i\leq N}
\mapsto\brk1{\min\big\{\max\{\xi_i,0\},255\big\}}_{1\leq i\leq N}.
\end{equation}
Moreover, \cite[Example~24.28]{Livre1} yields
\begin{equation}
\prox_{\gamma g_1}\colon\RR^N\to\RR^N\colon x\mapsto
\begin{cases}
x+\dfrac{\gamma}{d_{E}(x)}\brk1{\proj_Ex-x},&\text{if}\,\,
d_{E}(x)>\gamma;\\
\proj_{E}x,&\text{if}\,\,d_{E}(x)\leq\gamma,
\end{cases}
\end{equation}
where
\begin{equation}
\label{e:projE}
\proj_E x=\idft(\widetilde{x}),
\quad\text{with}\quad
\widetilde{x}(\nu)=
\begin{cases}
r(\nu), & \text{if}\;\nu\in R;\\
\widehat{x}(\nu), & \text{otherwise}.
\end{cases}
\end{equation}
Further, \cite[Proposition~24.11]{Livre1} yields
\begin{equation}
\prox_{\gamma g_2}\colon\RR^N\times\RR^N\to\RR^N\times\RR^N\colon
(\xi_i,\eta_i)_{1\leq i\leq N}\mapsto
\brk1{\varrho_i\xi_i,\varrho_i\eta_i}_{1\leq i\leq N},\,\,
\text{where}\,\,
\varrho_i=1-\frac{\sqrt{8}\gamma}{
\max\bigl\{\sqrt{8}\gamma,\|(\xi_i,\eta_i)\|\bigr\}}.
\end{equation}
At last, \cite[Example~2.3]{Siim19} establishes that
\begin{equation}
\nabla h\colon\RR^N\to\RR^N\colon x\mapsto
\dfrac{\rho_1H_1^*(H_1x-z_1)}{\max\{\rho_1,\|H_1x-z_1\|\}}+
\dfrac{\rho_2H_2^*(H_2x-z_2)}{\max\{\rho_2,\|H_2x-z_2\|\}}.
\end{equation}

\begin{figure}[t]
\centering
\setlength{\tabcolsep}{0.1pt} 
\begin{tabular}{ccc}
\includegraphics[width=4.6cm,height=4.5cm,clip]{
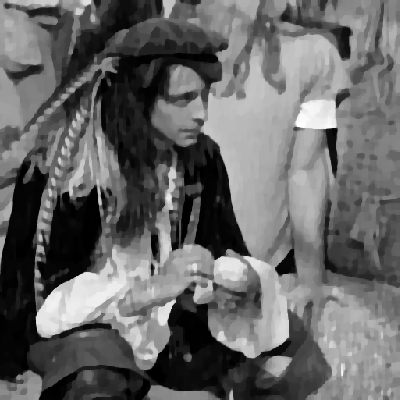}
&\includegraphics[width=4.6cm,height=4.5cm,clip]{
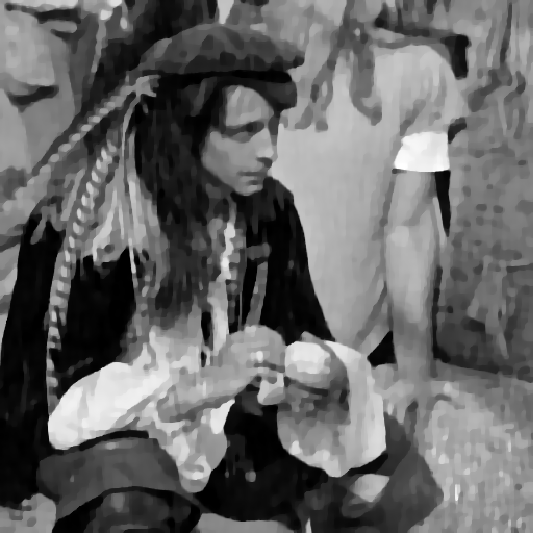}
&\includegraphics[width=4.6cm,height=4.5cm,clip]{
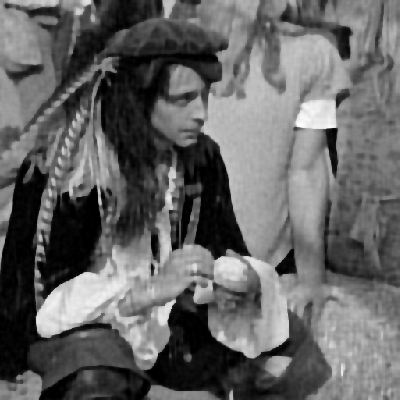}
\\
\small{(a) Problem~\ref{p:ex3a}.} 
&\small{(b) Problem~\ref{p:ex3c} ($\gamma=0.1$).} 
&\small{(c) Problem~\ref{p:ex3c} ($\gamma=0.99$).} 
\\ [0.3em]
\end{tabular} 
\caption{
Images reconstructed by Problems~\ref{p:ex3a} and \ref{p:ex3c}.}
\label{fig:ex3im2}
\end{figure}

\begin{figure}
\begin{tikzpicture}[scale=0.6]
\begin{axis}[height=8cm,width=23cm,legend style={at={(0.98,0.98)},
legend cell align={left},
anchor=north east}, xmin =0, xmax=90, xtick distance = 10, 
ytick distance = 10, ymin=-60, ymax=0, grid=both]
\addplot [line width=0.6mm, mark=none, color=dgreen] 
table[x={time}, y={error}]{figures/ex3/results/err_pd3o.txt};
\addlegendentry{Problem~\ref{p:ex3a}/Proposition~\ref{p:ave3}}
\addplot [line width=0.6mm, mark=none, color=magenta]table[
x={time},y={error}]{figures/ex3/results/gamma099/err_pcm099.txt};
\addlegendentry{Problem~\ref{p:ex3c}/Proposition~\ref{p:pcm}
($\gamma=0.99$)}
\end{axis}
\hfill
\end{tikzpicture}
\centering
\caption{Normalized error
$20\log_{10}(\|x_n-x_\infty\|/\|x_0-x_\infty\|)$ (dB) versus
time (s).
}
\label{fig:ex3im3}
\end{figure}

We construct the solution to Problem~\ref{p:ex3a} shown in 
Fig.~\ref{fig:ex3im2}(a) via Proposition~\ref{p:ave3}, where 
$\eta=1.99\beta$ and $\sigma=0.99/\eta$.
To apply Proposition~\ref{p:pcm} to Problem~\ref{p:ex3c}, let us
verify that condition \eqref{e:ppcm} is satisfied. 
By Proposition~\ref{p:1}\ref{p:1iiib} and
\cite[Corollary~16.48(iii)]{Livre1},
\begin{equation}
\partial\brk2{f+
\pcmx{\gamma}\brk1{d_E,\Id;\sqrt{8}\|\cdot\|_{1,2},D/\sqrt{8}}+h}
=\partial f+\partial\,
\pcmx{\gamma}\brk1{d_E,\Id;\sqrt{8}\|\cdot\|_{1,2},D/\sqrt{8}}
+\nabla h,
\end{equation} 
which is a maximally monotone operator with domain $C$
\cite[Theorem~20.25]{Livre1}. Hence, \eqref{e:ppcm} follows from
\cite[Corollary~21.25]{Livre1}. The image reconstructed by
\eqref{e:13} for Problem~\ref{p:ex3c} with $\gamma=0.1$ and 
$\lambda_n\equiv 1.89<1.90=2-\gamma/(2\beta)$ is shown in
Fig.~\ref{fig:ex3im2}(b). As predicted by
Theorem~\ref{t:4}\ref{t:4ivb}, since $\gamma$ is small, this
solution is similar to that produced by Problem~\ref{p:ex3a}
in Fig.~\ref{fig:ex3im2}(a). The solution
produced by Problem~\ref{p:ex3c} for $\gamma=0.99$ with 
$\lambda_n\equiv 1<1.01=2-\gamma/(2\beta)$ in \eqref{e:13} is shown
in Fig.~\ref{fig:ex3im2}(c) to yield a slightly sharper
reconstruction. Finally, Fig.~\ref{fig:ex3im3} illustrates the
faster convergence of algorithm~\eqref{e:13} for the proximal
comixture model \eqref{e:ex3c}.

\subsection{Experiment 3: Image reconstruction from phase}
\label{sec:ex5}

\begin{figure}[t]
\centering
\setlength{\tabcolsep}{0.2pt} 
\noindent
\begin{tabular}{cc}
\includegraphics[width=4.6cm,height=4.5cm,clip]{
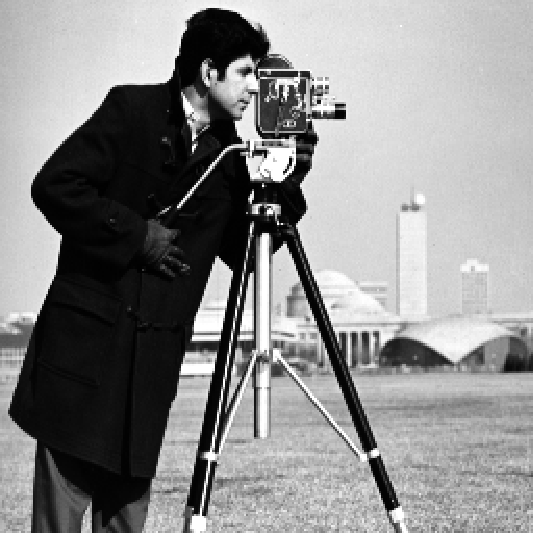} 
&\includegraphics[width=4.6cm,height=4.5cm,clip]{
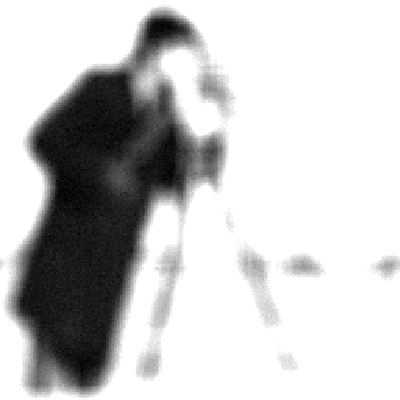}
\\
\small{(a)} & \small{(b)}
\\[0.3cm]
\includegraphics[width=4.6cm,height=4.5cm,clip]{
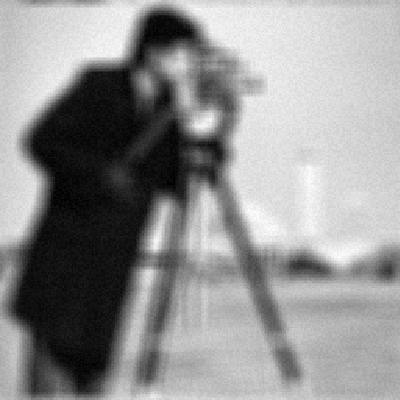} 
&\includegraphics[width=4.6cm,height=4.5cm,clip]{
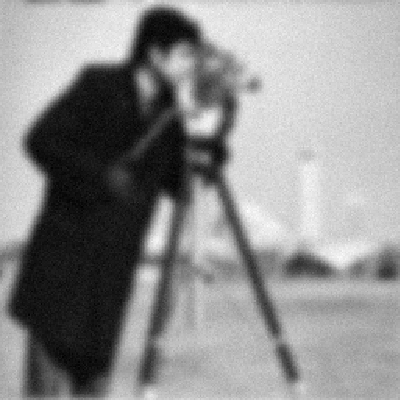} 
\\
\small{(c)} & \small{(d)} 
\end{tabular} 
\vskip -0.2cm
\caption{
(a) Original image $\overline{x}$. 
(b) Reference image $r$.
(c) Degraded image $z_1$.
(d) Degraded image $z_2$.}
\label{fig:ex5im1}
\end{figure}

We address a phase recovery problem considered in \cite{Siim19}.
The goal is to recover the original image 
$\overline{x}\in C=[0,255]^N$ ($N=256^2$) shown in
Fig.~\ref{fig:ex5im1}(a) from an imprecise observation of its
Fourier phase $\theta\approx\angle\,\widehat{\overline{x}}$
\cite{Levi83}. The problem is modeled as a convex feasibility
problem with the following constraint sets.
\begin{itemize}
\setlength{\itemsep}{-1pt}
\item
Phase: 
$C_1=\menge{x\in\RR^N}{\angle\;\widehat{x}=\theta}$.
\item
Proximity to the reference image $r$ of Fig.~\ref{fig:ex5im1}(b):
$C_2=\menge{x\in\RR^N}{\|x-r\|_2\leq\xi}$.
The image $r$ is a blurred and noise-corrupted version of
$\overline{x}$, which is further degraded by saturation (the pixel
values beyond 130 are clipped to 130).
\item
Upper bound on the norm of the gradient: 
$Dx/\sqrt{8}\in C_3$, where $D$ is defined as in \eqref{e:dgrad}
and $C_3=\menge{y\in\RR^N\times\RR^N}{\|y\|_2\leq\rho}$.
\item
A blurred observation $z_1=H_1\overline{x}+w_1$ is available 
(see Fig.~\ref{fig:ex5im1}(c)), where $H_1$ models a convolutional
blur with a uniform rectangular kernel of size $11\times 13$ and
$w_1$ is a realization of a bounded random noise vector with
components in $[-\eta,\eta]$, where $\eta=7$. This yields a 
blurred image-to-noise ratio of $31.90$ dB. A candidate solution
$x\in\RR^N$ should be such that, for every $k\in\{1,\ldots,N\}$,
the $k$th component of $H_1x-z_1$ lies in $[-\eta,\eta]$, that is
\begin{equation}
\scal{H_1x}{e_k}\in 
C_{k+3}=\bigl[\scal{z_1}{e_k}-\eta,\scal{z_1}{e_k}+\eta\bigr],
\end{equation}
where $(e_k)_{1\leq k\leq N}$ is the standard orthonormal basis of
$\RR^N$.
\item
A second blurred observation of $\overline{x}$ is available,
namely (see Fig.~\ref{fig:ex5im1}(d)) $z_2=H_2\overline{x}+w_2$,
where $H_2$ models a $17\times 17$ convolutional blur 
obtained by truncating the support of a Gaussian kernel 
with standard deviation $3$,
and $w_2$ is a realization of an unknown white 
Gaussian noise. The blurred image-to-noise ratio is $30.12$ dB.
This information is penalized by the function
$h\colon x\mapsto\norm{H_2x-z_2}^2/2$.
\end{itemize}
Because of inaccuracies in the values $\theta$, $\xi$, and
$\rho$, the convex feasibility problem arising from the above
constraints might be inconsistent and we relax it using the Berhu
function given by 
\begin{equation}
\mathfrak{b}\colon\RR\to\RR\colon\xi\mapsto
\begin{cases}
\dfrac{\xi^2+1}{2},&\text{if}\,\,\lvert\xi\rvert>1;\\
\lvert\xi\rvert,&\text{if}\,\,\lvert\xi\rvert\leq 1.
\end{cases}
\end{equation}
In addition, we set $L_1=L_2=\Id_{\RR^N}$, $L_3=D/\sqrt{8}$, and
$\alpha_1=\alpha_2=\alpha_3=1/4$ while, for every
$k\in\{1,\ldots,N\}$, 
$L_{k+3}\colon\RR^N\to\RR\colon x\mapsto\scal{H_1x}{e_k}$ and
$\alpha_{k+3}=1/(4N)$.

Our first formulation employs a standard composite average.

\begin{problem}
\label{p:ex5a}
The task is to
\begin{equation}
\label{e:ex5a}
\minimize{x\in C}{
\sum_{k=1}^{N+3}\alpha_k\mathfrak{b}\brk1{d_{C_k}(L_kx)}
+\dfrac{1}{2}\|H_2x-z_2\|^2}.
\end{equation}
\end{problem}

Our second formulation employs a proximal comixture.

\begin{problem}
\label{p:ex5c}
Let $\gamma\in\intv[o]{0}{2}$. The task is to
\begin{equation}
\label{e:ex5c}
\minimize{x\in C}{
\brk2{\pcmx{\gamma}\brk1{\mathfrak{b}\circ d_{C_k},L_k}_{1\leq
k\leq N+3}}(x)
+\dfrac{1}{2}\|H_2x-z_2\|^2}.
\end{equation}
\end{problem}

\begin{figure}[t]
\centering
\setlength{\tabcolsep}{0.1pt} 
\begin{tabular}{ccc}
\includegraphics[width=4.6cm,height=4.5cm]{
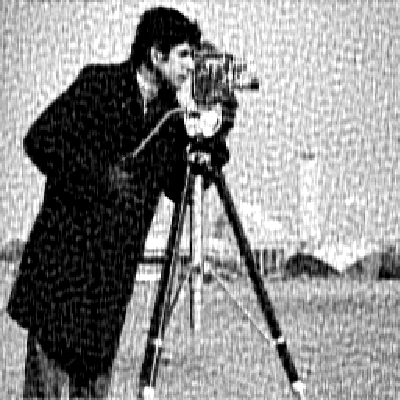}
&\includegraphics[width=4.6cm,height=4.5cm,clip]{
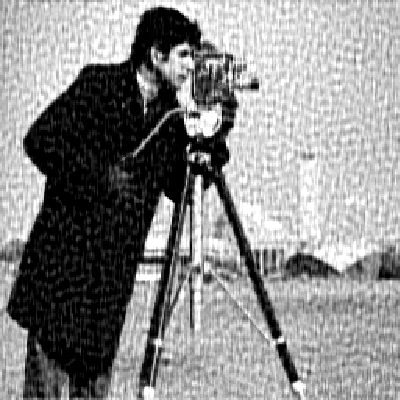}
&\includegraphics[width=4.6cm,height=4.5cm,clip]{
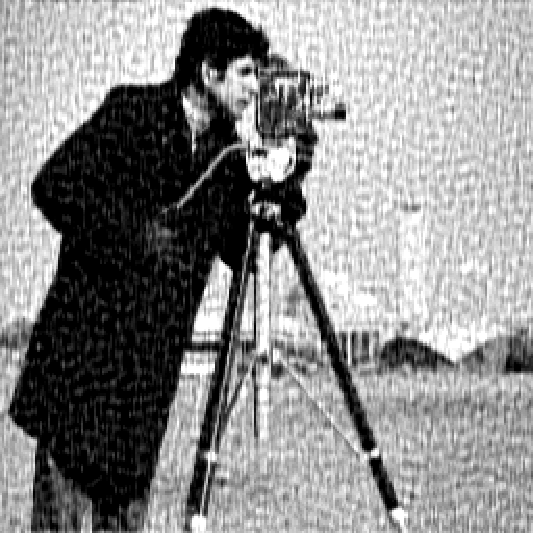}
\\
\small{(a) Problem~\ref{p:ex5a}.} 
&\small{(b) Problem~\ref{p:ex5c} ($\gamma=0.1$).} 
&\small{(c) Problem~\ref{p:ex5c} ($\gamma=1.99$).} 
\\ [0.3em]
\end{tabular} 
\caption{
Images reconstructed by Problems~\ref{p:ex5a} and
\ref{p:ex5c}.}
\label{fig:ex5im2}
\end{figure}

\begin{figure}
\begin{tikzpicture}[scale=0.6]
\begin{axis}[height=8cm,width=23cm,legend style={at={(0.98,0.98)},
legend cell align={left},
anchor=north east}, xmin =0, xmax=11000, xtick distance = 1000,
scaled x ticks=false,
ytick distance = 10, ymin=-60, ymax=0, grid=both]
\addplot [line width=0.6mm, mark=none, color=dgreen] 
table[x={time}, y={error}]{figures/ex4/results/err_pd3o.txt};
\addlegendentry{Problem~\ref{p:ex5a}/Proposition~\ref{p:ave3}}
\addplot [line width=0.6mm, mark=none, color=magenta]table[
x={time},y={error}]{figures/ex4/results/gamma199/err_pcm199.txt};
\addlegendentry{Problem~\ref{p:ex5c}/Proposition~\ref{p:pcm}
($\gamma=1.99$)}
\end{axis}
\hfill
\end{tikzpicture}
\centering
\caption{Normalized error
$20\log_{10}(\|x_n-x_\infty\|/\|x_0-x_\infty\|)$ (dB) versus
time (s).
}
\label{fig:ex5im3}
\end{figure}

Problems~\ref{p:ex5a} and \ref{p:ex5c} are particular instances of
Problems~\ref{prob:1} and \ref{prob:3}, respectively, where
$\HH=\RR^N$, $f=\iota_C$, $h=\|H_2\cdot-z_2\|^2/2$, $\beta=1$,
$p=N+3$, and, for every $k\in\{1,\ldots,p\}$, 
$g_k=\mathfrak{b}\circ d_{C_k}$.
We have $(\forall k\in\{1,\ldots,p\})$ $\|L_k\|^2\leq 1$. Further, 
$\prox_{\gamma f}$ is given as in \eqref{e:pbox},
\cite[Proposition~24.27]{Livre1} yields
\begin{equation}
\label{e:proxg12}
(\forall k\in\{1,2\})\quad
\prox_{\gamma g_k}\colon\RR^N\to\RR^N\colon x\mapsto
\begin{cases}
\proj_{C_k}x+\dfrac{1}{1+\gamma}\brk1{x-\proj_{C_k}x},
&\text{if}\,\,d_{C_k}(x)>1+\gamma;\\
x+\dfrac{\gamma}{d_{C_k}(x)}\brk1{\proj_{C_k}x-x},&\text{if}\,\,
\gamma<d_{C_k}(x)\leq 1+\gamma;\\
\proj_{C_k}x,&\text{if}\,\,d_{C_k}(x)\leq\gamma,
\end{cases}
\end{equation}
as well as
\begin{equation}
\label{e:proxg3}
\prox_{\gamma g_3}\colon\RR^N\times\RR^N\to\RR^N\times\RR^N\colon 
y\mapsto
\begin{cases}
\proj_{C_3}y+\dfrac{1}{1+\gamma}\brk1{y-\proj_{C_3}y},
&\text{if}\,\,d_{C_3}(y)>1+\gamma;\\
y+\dfrac{\gamma}{d_{C_3}(y)}\brk1{\proj_{C_3}y-y},&\text{if}\,\,
\gamma<d_{C_3}(y)\leq 1+\gamma;\\
\proj_{C_3}y,&\text{if}\,\,d_{C_3}(y)\leq\gamma
\end{cases}
\end{equation}
and
\begin{equation}
(\forall k\in\{4,\ldots,N+3\})\,\,
\prox_{\gamma g_k}\colon\RR\to\RR\colon\xi\mapsto
\begin{cases}
\proj_{C_k}\xi+\dfrac{1}{1+\gamma}\brk1{\xi-\proj_{C_k}\xi},
&\text{if}\,\,d_{C_k}(\xi)>1+\gamma;\\
\xi+\dfrac{\gamma}{d_{C_k}(\xi)}\brk1{\proj_{C_k}\xi-\xi},
&\text{if}\,\,\gamma<d_{C_k}(\xi)\leq 1+\gamma;\\
\proj_{C_k}\xi,&\text{if}\,\,d_{C_k}(\xi)\leq\gamma,
\end{cases}
\end{equation}
where $\proj_{C_1}$ is given by \cite{Levi83}
\begin{equation}
\label{e:projC1}
\proj_{C_1}x=\idft\brk2{
\lvert\dft x\rvert\max\{\cos\brk1{\angle(\dft x)-\theta},0\}
\exp(\imath\theta)},
\end{equation}
and, for every $k\in\{2,\ldots,N+3\}$, $\proj_{C_k}$ follows from
\cite[Example~3.18 and Proposition~3.19]{Livre1}.
Finally, $\nabla h=H_2^*\circ(H_2\cdot-z_2)$.

We apply Proposition~\ref{p:ave3} to Problem~\ref{p:ex5a} with
$\eta=1.99\beta$ and $\sigma=0.99/\eta$, which produces the
reconstructed image shown in Fig.~\ref{fig:ex5im2}(a). Following
the same argument used in Problem~\ref{p:ex3c}, we note that
condition \eqref{e:ppcm} is satisfied and we apply
Proposition~\ref{p:pcm} to Problem~\ref{p:ex5c}. The image
reconstructed by \eqref{e:13} for $\gamma=0.1$ and
$\lambda_n\equiv 1.94<1.95=2-\gamma/(2\beta)$ is shown in
Fig.~\ref{fig:ex5im2}(b). This
solution is similar to that of Fig.~\ref{fig:ex5im2}(a), which is
consistent with Theorem~\ref{t:4}\ref{t:4ivb} since $\gamma$ is
small. For Problem~\ref{p:ex5c} with 
$\gamma=1.99$, algorithm
\eqref{e:13} with $\lambda_n\equiv 1<1.005=2-\gamma/(2\beta)$
yields the somewhat sharper reconstruction shown in
Fig.~\ref{fig:ex5im2}(c). Fig.~\ref{fig:ex5im3} depicts the faster
convergence of algorithm~\eqref{e:13} for the proximal comixture
model \eqref{e:ex5c}.

\subsection{Experiment 4: Linear regression}
\label{sec:ex1}

This experiment focuses on a linear regression model with an
overlapping group structure on the inputs \cite{Chen12,Yuan11}.
We consider $p=100$ groups of indices 
$(I_k)_{1\leq k\leq p}$ in $\{1,\ldots,N\}$ ($N=90p+10=9010$) of
length 100 such that two consecutive groups overlap by 10
variables, i.e.,
\begin{equation}
I_1=\{1,\ldots,100\},
I_{2}=\{91,\ldots,190\},\ldots,
I_{p}=\{N-99,\ldots,N\}.
\end{equation}
We use $M=10000$ samples.
The entries of the input matrix $A\in\RR^{M\times N}$ are i.i.d.
samples from a $\mathcal{N}(0,1)$ distribution. The entries of the
true regression coefficients $\overline{x}\in\RR^N$ are also i.i.d.
samples from a $\mathcal{N}(0,1)$ distribution, and the output data
is generated by the noisy linear model $z=A\overline{x}+w$, where
$w\in\RR^M$ has entries that are i.i.d. samples from a
$\mathcal{N}(0,1)$ distribution. For every $k\in\{1,\ldots,p\}$, we
set $L_k\colon\RR^N\to\RR^{100}\colon x
=(\xi_i)_{1\leq i\leq N}\mapsto(\xi_i)_{i\in I_k}$.

As before, we consider formulations based on a standard composite 
average and the proximal comixture.

\begin{problem}
\label{p:ex1a}
The task is to 
\begin{equation}
\label{e:ex1a}
\minimize{x\in\RR^N}{\dfrac{1}{p}\|x\|_1+
\dfrac{1}{p}\sum_{k=1}^{p}\|L_kx\|+\dfrac{1}{2p^2}\|Ax-z\|^2}.
\end{equation}
\end{problem}

\begin{problem}
\label{p:ex1c}
Let $\gamma\in\intv[o]{0}{2p^2/\|A\|^2}$. The task is to
\begin{equation}
\label{e:ex1c}
\minimize{x\in\RR^N}{\dfrac{1}{p}\|x\|_1+
\brk2{\pcmx{\gamma}\brk1{\|\cdot\|,L_k}_{1\leq k\leq p}}(x)
+\dfrac{1}{2p^2}\|Ax-z\|^2}.
\end{equation}
\end{problem}

Problems~\ref{p:ex1a} and \ref{p:ex1c} are particular instances of
Problems~\ref{prob:1} and \ref{prob:3}, respectively, where
$\HH=\RR^N$, $f=\|\cdot\|_1/p$, $h=\|A\cdot-z\|^2/(2p^2)$,
$\beta=p^2/\|A\|^2$, and, for every $k\in\{1,\ldots,p\}$,
$\alpha_k=1/p$, $\GG_k=\RR^{100}$, $g_k=\|\cdot\|$, and 
$\|L_k\|=1$. Additionally, 
\begin{equation}
\prox_{\gamma f}\colon\RR^N\to\RR^N\colon(\xi_i)_{1\leq i\leq N}
\mapsto\brk1{\sign{(\xi_i)}\max\{\lvert\xi_i\rvert-\gamma/p,0\}}_{
1\leq i\leq N}.
\end{equation}
Further, \cite[Example~24.20]{Livre1} yields
\begin{equation}
(\forall k\in\{1,\ldots,p\})\,\,
\prox_{\gamma g_k}\colon\RR^{100}\to\RR^{100}\colon y\mapsto
\brk3{1-\dfrac{\gamma}{\max\{\|y\|,\gamma\}}}y,
\end{equation}
while $\nabla h=(1/p^2)A^*\circ(A\cdot-z)$.
Next, let us verify that condition \eqref{e:ppcm} is satisfied. 
Note that, by Proposition~\ref{p:2}, 
$\pcmx{\gamma}(\|\cdot\|,L_k)_{1\leq k\leq p}\geq 0$.
Thus, the objective function of Problem~\ref{p:ex1c} is coercive
and the existence of minimizers is guaranteed by
\cite[Proposition~11.15(i)]{Livre1}. Therefore, Fermat's rule
\cite[Theorem~16.3]{Livre1} and \cite[Corollary~16.48(iii)]{Livre1}
guarantee that condition \eqref{e:ppcm} holds. 
The solutions $x_{\infty}$ to Problem~\ref{p:ex1a} and
to Problem~\ref{p:ex1c} for $\gamma=1.99\beta$ are very
close, with a relative difference 
$\|x_{\infty,\eqref{e:ex1a}}-x_{\infty,\eqref{e:ex1c}}\|/
\|x_{\infty,\eqref{e:ex1a}}\|\approx 5.045\times 10^{-6}$.
In addition, 
$\|x_\infty-\overline{x}\|/\|\overline{x}\|\approx 0.153$. 
For Problem~\ref{p:ex1a}, we use algorithm~\eqref{e:12} with
$\eta=1.99\beta$ and $\sigma=0.99/\eta$, while for
Problem~\ref{p:ex1c} for $\gamma=1.99\beta$, we use
algorithm~\eqref{e:13} with $\lambda_n\equiv
1<1.005=2-\gamma/(2\beta)$. The convergence pattern of algorithm
\eqref{e:13} for the proximal comixture model \eqref{e:ex1c} is
shown in Fig.~\ref{fig:ex1im1}.

\begin{figure}[ht!]
\begin{tikzpicture}[scale=0.6]
\begin{axis}[height=8cm,width=23cm,legend style={at={(0.98,0.98)},
legend cell align={left},
anchor=north east}, xmin =0, xmax=24, xtick distance = 3, 
ytick distance = 10, ymin=-60, ymax=0, grid=both]
\addplot [line width=0.6mm, mark=none, color=dgreen] 
table[x={time}, y={error}]{figures/ex1/err_pd3o.txt};
\addlegendentry{Problem~\ref{p:ex1a}/Proposition~\ref{p:ave3}}
\addplot [line width=0.6mm, mark=none, color=magenta] 
table[x={time}, y={error}]{figures/ex1/err_pcm.txt};
\addlegendentry{Problem~\ref{p:ex1c}/Proposition~\ref{p:pcm} 
($\gamma=1.99\beta$)}
\end{axis}
\hfill
\end{tikzpicture}
\centering
\caption{Normalized error
$20\log_{10}(\|x_n-x_\infty\|/\|x_0-x_\infty\|)$ (dB) versus
time (s).
}
\label{fig:ex1im1}
\end{figure}

\end{document}